\documentclass[11pt]{article}
\usepackage[dvips]{color}
\usepackage{mathrsfs}
\usepackage{tablefootnote}
\usepackage{graphicx}
\usepackage{amsfonts}
\usepackage{amsmath}
\usepackage{amssymb}
 \usepackage{epstopdf}
\definecolor{dgreen}{rgb}{0,.8,.3}
\definecolor{blue}{rgb}{.2,.3,.7}
\definecolor{red}{rgb}{1.0,0.2,0.2}
\usepackage{color}
\usepackage{pgfplots}
\usepackage{subcaption}
\usepackage{longtable,booktabs,easyReview}
\usepackage{multirow}
\input amssym.def

\input epsf.tex

\begin{document}

\renewcommand{\Box}{\rule{2.2mm}{2.2mm}}
\newcommand{\BOX}{\hfill \Box}

\newtheorem{eg}{Example}[section]
\newtheorem{thm}{Theorem}[section]
\newtheorem{lemma}{Lemma}[section]
\newtheorem{example}{Example}[section]
\newtheorem{remark}{Remark}[section]
\newtheorem{proposition}{Proposition}[section]
\newtheorem{corollary}{Corollary}[section]
\newtheorem{defn}{Definition}[section]
\newtheorem{alg}{Algorithm}[section]
\newtheorem{ass}{Assumption}[section]
\newenvironment{case}
    {\left\{\def\arraystretch{1.2}\hskip-\arraycolsep \array{l@{\quad}l}}
    {\endarray\hskip-\arraycolsep\right.}

\def\argmin{\mathop{\rm argmin}}

\makeatletter
\renewcommand{\theequation}{\thesection.\arabic{equation}}
\@addtoreset{equation}{section} \makeatother

\title{Enhanced Barrier-Smoothing Technique for Bilevel Optimization with Nonsmooth Mappings}
\author{Mengwei Xu\thanks{\baselineskip 9pt Institute of Mathematics, Hebei University of Technology, Tianjin 300401, China. 
E-mail: xumengw@hotmail.com. The
research of this author is supported by the National Natural
Science Foundation of China under Project No. 12071342 and the Natural Science Foundation of Hebei Province, No. A2020202030},\  
Yu-Hong Dai\thanks{\baselineskip 9pt Academy of Mathematics and Systems Science, Chinese Academy of Sciences, Beijing 100190, China \& School
of Mathematical Sciences, University of Chinese Academy of Sciences, Beijing 100049, China. The
research of this author is
supported by the NSFC grants (Nos. 12021001, 11991021, 11991020 and 11971372), the Strategic Priority Research
Program of Chinese Academy of Sciences (No. XDA27000000) and Beijing Academy of Artificial Intelligence},\ 
Xin-Wei Liu\thanks{\baselineskip 9pt Institute of Mathematics, Hebei University of Technology, Tianjin 300401, China. E-mail: mathlxw@hebut.edu.cn. The
research of this author is supported by the NSFC grants (No. 12071108 and No. 11671116)}\ 
and Bo Wang\thanks{\baselineskip 9pt Key Laboratory of Operations Research and Control of Universities in Fujian, School of Mathematics and Statistics, Fuzhou University, Fuzhou 350116, China. Email: bowang@fzu.edu.cn.
The
research of this author is supported by Fujian Provincial Natural Science Foundation of China (No. 2023J01416)
}
}
\date{}
\maketitle
{\bf Abstract.} Bilevel optimization problems, encountered in fields such as economics, engineering, and machine learning, pose significant computational challenges due to their hierarchical structure and constraints at both upper and lower levels. Traditional gradient-based methods are effective for unconstrained bilevel programs with unique lower level solutions, but struggle with constrained bilevel problems due to the nonsmoothness of lower level solution mappings. To overcome these challenges, this paper introduces the Enhanced Barrier-Smoothing Algorithm (EBSA), a novel approach that integrates gradient-based techniques with an augmented Lagrangian framework.
 EBSA utilizes innovative smoothing functions to approximate the primal-dual solution mapping of the lower level problem, and then transforms the bilevel problem into a sequence of smooth single-level problems.
 This approach not only addresses the nonsmoothness but also enhances convergence properties. 
Theoretical analysis demonstrates its superiority in achieving Clarke and, under certain conditions, Bouligand stationary points for bilevel problems.
Both theoretical analysis and preliminary numerical experiments confirm the robustness and efficiency of EBSA.



\newpage

\baselineskip 18pt
\parskip 2pt
\section{Introduction}
Bilevel programming problems are a specialized class of optimization problems characterized by a hierarchical structure, where the feasible region of the upper level problem is contingent on the solution set of the lower level problem. This unique structure has led to widespread applications in domains such as Stackelberg games, principal-agent problems \cite{Mirrlees99, s}, hyperparameter optimizations \cite{bkhp, lmyzz, lmyzz2, lvd}, meta-learning, and neural architecture search \cite{ffsgp, rfkl}. 
Despite their broad applicability, bilevel programs pose significant computational challenges due to their intrinsic complexity, nonconvexity and  nondifferentiable nature.

The mathematical formulation of a bilevel optimization problem is as follows:
\begin{equation}
	\begin{array}{ll}\displaystyle
		\min_{x,y} & F(x,y)\\
		{\rm s.t.}&G(x,y)\leq 0,\ H(x,y)= 0,\\
		&y\in S(x)
	\end{array}
\tag{BP}
\end{equation}
where  $S(x)$ denotes the solution set of the lower level program
	\begin{equation}
		\begin{array}{ll}\displaystyle
			\min_{y} & f(x,y)\\
			{\rm s.t.}&g(x,y)\leq 0.
		\end{array}
		\tag{P$_{x}$}
	\end{equation}
Here we assume that 
$F:\mathbb{R}^d\times \mathbb{R}^l \rightarrow \mathbb{R}$ 
, $G:\mathbb{R}^d\times \mathbb{R}^l \rightarrow \mathbb{R}^p$ and $H:\mathbb{R}^d\times \mathbb{R}^l \rightarrow \mathbb{R}^q$  
are continuously differentiable,  
$f:\mathbb{R}^d\times \mathbb{R}^l \rightarrow \mathbb{R}$ and $g:\mathbb{R}^d\times \mathbb{R}^l \rightarrow \mathbb{R}^m$
are continuously differentiable and twice continuously differentiable with respect to the variable $y$.

\subsection{Existing Solution Methods}
 Numerous methods have been developed to tackle bilevel optimization problems.
Gradient-based methods, in particular, have demonstrated significant efficiency in solving unconstrained bilevel programs. When problem (BP) is unconstrained and  $({\rm P}_x)$ has a unique solution for each $x$,
it can be reformulated as $\min_x\ \phi(x):=F(x,y(x))$, where $y(x)= \argmin_y\ f(x,y)$.
By utilizing the chain rule, we can express the gradient as:
\begin{eqnarray*}
\nabla \phi(x)=\nabla_x F(x,y(x)) + \nabla y(x)^T \nabla_y F(x,y(x)),
\end{eqnarray*}
where $\nabla y(x)$ can be obtained using the Implicit Function Theorem
 \cite{gfps, lmyzz} or
approximated through a dynamic system \cite{fdfp, ffsgp, gfps}.
However, applying these methods to constrained bilevel programs is challenging due to the lack of continuous differentiability of the lower level solution mapping \cite{dz}.
This nonsmoothness adversely affects the convergence and stability of the optimization process, making it difficult to achieve reliable and efficient solutions.
Recent study \cite{kpwzs} attempted to extend gradient-based methods to bilevel problems where the upper level has abstract constraints and the lower level has linear constraints, under the assumption that the lower level solution mapping is Lipschitz continuous.

The first-order approach has gained considerable attention by substituting the solution set of the lower level problem with its necessary optimality conditions, resulting in a formulation known as the mathematical program with equilibrium constraints (MPEC) \cite{as, lpr, okz}.
However, Dempe et al. \cite{DempeDu, df16} demonstrated that the MPEC reformulation and the bilevel problem are not equivalent in terms of local solutions, even when the lower level problem is convex with respect to the lower level variable. 

Building on the value function reformulation proposed by Outrata \cite{Outrata}, 
Lin et al. introduced smoothing methods 
for nonconvex bilevel program, 
 where the constraint set of $({\rm P}_x)$ is independent of the upper level variable \cite{lxy, xyz}.
Based on partial calmness, researchers \cite{fzz,jmz,tz} introduced  an equation system for the value function reformulation, parameterized by partial exact penalization. They developed algorithms such as the semi-smooth Newton and Levenberg–Marquardt methods to solve it.
The Moreau envelope reformulation and the use of difference of convex (DC) algorithms were proposed for bilevel problems where $({\rm P}_x)$ is convex \cite{gttzz, yyzz}.
However, solving the value function reformulation remains difficult due to its implicit nature, nonsmoothness, and the challenge of satisfying usual constraint qualifications at any feasible point. Consequently, dealing with the value function reformulation remains a formidable task.
Further numerical developments for bilevel programs with special structures can be found in \cite{dkpk,dz, hss22,hxlt,llz,lm,nwy,ss17} and the references therein.

\subsection{Research Motivation and Contributions}
Despite significant advancements, existing methods for solving bilevel problems face several critical challenges:
\begin{itemize}
\item[{\rm (i)}] Nondifferentiability and strict complementarity: Traditional gradient-based methods for solving constrained bilevel programs, struggled with the nondifferentiability in the lower level solution mapping.
This issue arises from the strict complementarity (SC) condition within the Karush-Kuhn-Tucker (KKT) system of $({\rm P}_x)$, which is too stringent and difficult to satisfy. 
 
 \item[{\rm (ii)}] Lack of optimality guarantees: First-order approaches (MPEC reformulation) are not fully equivalent to bilevel programs in terms of local solutions. The equivalence holds only under restrictive assumptions \cite{DempeDu}.
 
\item[{\rm (iii)}] Computational difficulties: The MPEC reformulation and the value function reformulation encounter considerable obstacles due to the failure of traditional constraint qualifications.  Consequently, they are inherently complex optimization problems, leading to a substantial increase in computational effort and time consumption.
\end{itemize}

To address these challenges, we propose the Enhanced Barrier-Smoothing Algorithm (EBSA) for solving constrained bilevel programs.
We propose a family of smoothing functions to approximate the lower level primal-dual solution mapping while preserving gradient consistency. These functions are derived from the  smoothing barrier augmented Lagrangian (SBAL) approach \cite{ld} (see Section 2.2 for more details). 
By combining the strengths of logarithmic barrier and augmented Lagrangian techniques, SBAL effectively handles inequality constraints, eliminating the need for primal and dual iterates to be interior points.
Specifically, we reformulate the bilevel problem as
\begin{equation}
		\begin{array}{ll}\displaystyle
			\min_{x} & F(x,y(x))\\
			{\rm s.t.}&G(x,y(x))\leq 0,\\
			&  H(x,y(x))= 0,
		\end{array}
		\tag{\rm SP}
	\end{equation}
where $y(x)$ is the solution function of the lower level problem.
This reformulation allows for a hybrid approach, combining the robustness of gradient-based methods with the flexibility of augmented Lagrangian methods.

Without assuming the SC condition, Dai and Zhang \cite{dz20} discovered that the primal-dual solution mapping demonstrates Lipschitz continuous differentiability under the linear independence constraint qualification (LICQ) and the strong second-order sufficient optimality condition (SSOSC). In addition, they provided a detailed characterization of the Bouligand subdifferentials (B-subdifferentials) and the Clarke subdifferentials (C-subdifferentials) of the primal-dual solution mapping, as shown in Lemma \ref{implicit-sub}.  We define the B-stationary points and C-stationary points of the problem (SP) accordingly.

Our contributions are as follows:
\begin{itemize}
\item[{\rm (i)}] Innovative smoothing functions: We introduce a family of smoothing functions derived from the SBAL method. These functions approximate the primal-dual solution mapping of the lower level problem while maintaining gradient consistency, thus enhancing the applicability of gradient-based methods to constrained bilevel problems.
\item[{\rm (ii)}] Enhanced robustness: We develop an innovative line search procedure within our method, significantly improving the robustness of the algorithm. 
\item[{\rm (iii)}] Convergence theory: Our theoretical analysis demonstrates that any accumulation point of the iteration sequence generated by EBSA corresponds to a C-stationary point or a B-stationary point of (SP), provided the multipliers are bounded. 
\end{itemize}

The paper is structured as follows. 
In Section 2, we provide an overview of the necessary preliminaries, including the differentiability results of the primal-dual solution mapping and the SBAL method for the lower level problem.
Section 3 details our proposed smoothing functions for the primal-dual solution mapping, along with the development of the gradient consistency property.
In Section 4, we introduce the EBSA method and analyze the convergence properties of the algorithm.
Section 5 showcases numerical experiments conducted on various bilevel problems to validate the effectiveness of our approach.
Finally, Section 6 concludes the paper.

We adopt the following standard notation in this paper. For any two vectors $a$ and $b $ in $\mathbb{R}^n$, we denote by either $\langle a, b \rangle$ or  $a^T b$  their inner product. Given a function $G: \mathbb{R}^n\rightarrow \mathbb{R}^m$, we denote its Jacobian by $\nabla G(z)\in \mathbb{R}^{m\times n}$ and, if $m=1$, the gradient $\nabla G(z)\in \mathbb{R}^n$ is considered as a column vector.   Denote by $I_{m}$ the $m\times m$ identity matrix.
For a set $\Omega\subseteq \mathbb{R}^n$, we denote by 
 co $\Omega$ the convex hull of $\Omega$.
For a matrix $A\in \mathbb{R}^{n\times m}$, $A^T$ denotes its transpose. 
Let $e_i\in \mathbb{R}^m$ be the vector such that the $i$th component is 1 and others are 0.

\section{Preliminaries}
In this section, we provide a brief introduction to the differentiability properties of the primal-dual solution mapping for the lower level problem, as well as the SBAL method for optimization problem with inequality constraints \cite{ld}.

\subsection{Differentiability of the primal-dual solution mapping}
For  a Lipschitz continuous function $\phi: \mathbb{R}^d \rightarrow  \mathbb{R}$ at $\bar  x$,
 the Bouligand subdifferential of $\phi $ at $ \bar  x$ is denoted by
 $\partial_B \phi(x) := \{ v \in \mathbb{R}^n \mid \exists x_k\in D_f, x_k \to x, \nabla f(x_k) \to v \text{ as } k \to \infty \}$, where $D_f$ denotes the set of points where $f$ is differentiable.
The  Clarke subdifferential of $\phi$ at ${x}$ is
$$ \partial_C \phi({x}):={co} \partial_B \phi({x}).$$
Detailed discussions on the subdifferentials can be found in \cite{var}. 


For the constrained problem (BP), consider any $x$, let $(y,u)$ be a KKT point of the lower level problem $({\rm P}_x)$. This gives:
\begin{eqnarray*}
\nabla_y\mathcal{L}(x,y,u)=0,\quad 
 \min\{u, -g(x,y)\}=0,
\end{eqnarray*} 
where $\mathcal{L}(x,y,u)=f(x,y)+  \sum_{i=1}^m u_i  g_i(x,y)$.
Denote the active index set by $I_x(y):=\{i\in\{1,\cdots,m\}: g_i(x,y)=0\}$.

Based on the Implicit Function Theorem, for any $x$, there exists unique solution $(y(x),u(x))$ for the  problem $({\rm P}_x)$, which is continuously differentiable  under the LICQ, the SC condition, i.e., $\bar{u}_i - g_i(\bar{x},\bar{y})>0,\ \forall i=1,\cdots,m$ and the second-order sufficient optimality condition for $P_x$,
see e.g. \cite[Theorem 2.1]{r74}.
Since the SC condition is too strong to be satisfied, Dai and Zhang \cite{dz20} studied the Lipschitz continuity of $(y(x),u(x))$ under the LICQ and the SSOSC.

\begin{ass}\label{strongregular}
Assume that $(\bar{y},\bar{u})$ satisfies the LICQ and the SSOSC  for problem $(P_{\bar x})$, i.e., 
\begin{eqnarray*}
\langle \nabla_{yy}^2 \mathcal{L}(\bar{x},\bar{y},\bar{u}) d_y, d_y  \rangle>0,\ \forall d_y\in {\rm aff}\ C_{\bar{x}}(\bar{y})\setminus \{0\},
\end{eqnarray*} 
where $C_x(y):=\{d_y\in\mathbb{R}^l: \nabla_y g_i(x,y)^T d_y\leq 0,\ i\in I_x(y);\ \nabla_y f(x,y)^T d_y\leq 0\}$ denotes the critical cone of $({\rm P}_x)$ at $y$.
\end{ass}

\begin{lemma}\cite[Lemma 2.2, Proposition 2.4]{dz20}\label{implicit-sub}
Assume that $(\bar{y},\bar{u})$ is a KKT point of $({\rm P}_{\bar{x}})$ satisfying the Assumption \ref{strongregular}.
Then there exist $\delta_0>0, \varepsilon_0>0$ and a locally Lipschitz continuous mapping $(y(x),u(x)):  \mathbb{B}_{\delta_0}(\bar{x})\to  \mathbb{B}_{\varepsilon_0}(\bar{y})\times \mathbb{B}_{\varepsilon_0}(\bar{u})$, which satisfies 
the KKT condition of $({\rm P}_{x})$ and the Assumption \ref{strongregular}. The  
Bouligand subdifferential (B-subdifferential) and the Clarke subdifferential (C-subdifferential) of $(y(x),u(x))$ can be respected as 
\begin{eqnarray*}
\partial_B (y(x), u(x))^T
  \subseteq
M^B(x):=
 \left \{ - \mathcal{A}(x,W)^{-1} a(x,W): W\in \partial_B \Pi_{\mathbb{R}_-^m}(u(x)+g(x,y(x)))
   \right\}
\end{eqnarray*}
and 
\begin{eqnarray*}
\partial_C (y(x), u(x))^T
  \subseteq
M(x):=
 \left \{ - \mathcal{A}(x,W)^{-1} a(x,W): W\in \partial_C \Pi_{\mathbb{R}_-^m}(u(x)+g(x,y(x)))
   \right\},
\end{eqnarray*}
respectively, where $\Pi_{\mathbb{R}_-^m}(u(x)+g(x,y(x)))$ is the projector of $u(x)+g(x,y(x))$ onto $\mathbb{R}_-^m$,
\begin{eqnarray*}
\mathcal{A}(x,W):=\left (\begin{array}{ll}
\displaystyle\nabla_{yy}^2 f(x,y(x))+  \sum_{i=1}^m u_i(x) \nabla_{yy}^2 g_i(x,y(x))& 
\nabla_{y} g(x,y(x))^T\\
~~~~~~~~~~~~ (W-I)\nabla_y g(x,y(x))& ~~~~~~W
  \end{array} \right)
\end{eqnarray*} 
and 
\begin{eqnarray*}
a(x,W):=\left (\begin{array}{l}
\nabla_{xy}^2 f(x,y(x))+  \sum_{i=1}^m u_i(x) \nabla_{xy}^2 g_i(x,y(x))\\
~~~~~~~~~ (W-I)\nabla_x g(x,y(x))
  \end{array} \right).
\end{eqnarray*} 
\end{lemma}  


The rest of this subsection reviews optimality conditions of the problem (SP). 
\begin{ass}\label{EMF-U}
We assume that the extended no nonzero abnormal multiplier constraint qualification $(\rm ENNAMCQ)$ holds at $ \bar x$ with $y(\bar x)=\bar y$ for $({\rm SP})$
if for $\lambda_i\geq 0,\  i=1,\cdots,p,$
\begin{eqnarray*}
&&0\in \nabla_x G(\bar x,\bar y)^T\lambda+ \nabla_x H(\bar x,\bar y)^T\mu+M_{y}(\bar x)^T [ \nabla_y G(\bar x,\bar y)^T \lambda +\nabla_y H(\bar x,\bar y)^T\mu],\\
&& \sum_{i=1}^p {\lambda}_i G_i(\bar x,\bar y)+\sum_{j=1}^q \mu_j H_j(\bar x,\bar y)\geq 0,
\end{eqnarray*}
imply that $ \lambda=0,\mu=0$. Here $M_{y}( x):=\{-[I_m,0] \mathcal{A}(x,W)^{-1} a(x,W):W\in \partial_C \Pi_{\mathbb{R}_-^m}(u(x)+g(x,y(x)))\}$ represents the upper bound of $\partial_C y(x)$.
\end{ass}
If $\bar{x}$ is feasible for $({\rm SP})$, the ENNAMCQ reduces to the NNAMCQ, which 
 equals to the nonsmooth version of Mangasarian Fromovitz constraint qualification (MFCQ).
Under the Assumptions \ref{strongregular}, \ref{EMF-U}, we derive the following optimality conditions.

\begin{thm}\label{gradF-r}\cite[Theoerm 5.4]{Liu 24}
Let $(\bar{x},\bar{y})$ be a local minimizer of (SP).
Assume that the Assumptions \ref{strongregular}, \ref{EMF-U} hold.
 Then $(\bar{x},\bar{y})$ is a C-stationary point of $({\rm SP})$, i.e., there exist $\lambda\geq 0,\mu$
 such that $\lambda_i G_i(\bar x,\bar y)=0$ for each $i=1,\cdots,p$ and
\begin{eqnarray*}
0&\in&\nabla_x F(\bar x,\bar y) + \nabla_x G(\bar x,\bar y)^T\lambda+ \nabla_x H(\bar x,\bar y)^T\mu\\
&&
+M_{y}(\bar x)^T [\nabla_y F(\bar x,\bar y)+ \nabla_y G(\bar x,\bar y)^T \lambda +\nabla_y H(\bar x,\bar y)^T\mu].
\end{eqnarray*}
If $M_{y}(\bar x)$ is replaced to $M_{y}^B(\bar x)$, where $M_y^B(x):=\{-[I_m,0] \mathcal{A}(x,W)^{-1} a(x,W):W\in \partial_B \Pi_{\mathbb{R}_-^m}(u(x)+g(x,y(x)))\}$  represents the upper bound of $\partial_B y(x)$, we say
 $(\bar{x},\bar{y})$ is a B-stationary point.
\end{thm}  

\subsection {The SBAL method  for the constrained lower level problem}
In this subsection, we apply the SBAL method \cite{ld,ldhs} to the lower level problem $({\rm P}_x)$. 
For $r>0$, by introducing auxiliary variables $z_i> 0, i=1,\cdots,m$, $({\rm P}_x)$ can be approximated by 
\begin{equation}
	\begin{array}{ll}\displaystyle
		\min_{y,z} &  f(x,y)-r \sum_{i=1}^m \ln z_i\\
		{\rm s.t.} & z_i+g_i(x,y)=0, i=1,\cdots,m.
	\end{array}
	\tag{${\rm P}_x^{r}$}
\end{equation}
For each $\rho>0$, $s\in\mathbb{R}^m$, 
the augmented Lagrangian function of $({\rm P}_x^{r})$ is defined as follows:
\begin{eqnarray*}
f_r^{\rho}(x,y,z,s):= f(x,y)+ \sum_{i=1}^m [-r \ln z_i+ s_i(z_i+g_i(x,y))+\frac{1}{2 \rho}(z_i+g_i(x,y))^2].
\end{eqnarray*}
To ensure that $s$ is a good estimate of Lagrange multiplier vector, the authors in \cite{ld,ldhs} maximized the augmented Lagrangian with respect to $s$. This leads to the following problem:
\begin{equation}\label{minmax}
	\tag{${\rm LP}_x^{r,\rho}$}
	\min_{y,z} \max_s ~~  f_r^{\rho}(x,y,z,s).
\end{equation}
For each $x$ and any solution $(y,z,s)$ of the problem $({\rm LP}_x^{r,\rho})$, we must have that $\nabla_z f_r^{\rho}(x,y,z,s)=0$, 
which derives that $z_i$ is a function depending on $x,y,s_i,r,\rho$\alert{:}
 \begin{eqnarray}\label{z}
z_i(x,y,s_i,r,\rho):= \frac{1}{2}[\sqrt{(\rho s_i+ g_i(x,y))^2+4 r \rho}-(\rho s_i+g_i(x,y))],\ i=1,\cdots,m.
\end{eqnarray}
Since the function $f_r^{\rho}$ is convex with respect to the variable $z$, we can replace $z$ by $z(x,y,s,r,\rho)$ 
and rewrite $f_r^{\rho}(x,y,z,s)$ by $\tilde{f}_r^{\rho}(x,y,s)$. 

For any $x$ and $r>0, \rho>0$,
any solution of $({\rm LP}_x^{r,\rho})$ satisfies the KKT conditions, 
\begin{eqnarray*}
&&\phi^{r,\rho}(x,y,s):=\nabla_y \tilde{f}_r^{\rho}(x,y,s)=\nabla_{y} f(x,y)+  \sum_{i=1}^m \frac{\kappa_i}{\rho} \nabla_{y} g_i(x,y) =0,\\
&&\psi^{r,\rho}(x,y,s):=\nabla_s \tilde{f}_r^{\rho}(x,y,s)=z(x,y,s,r,\rho)+g(x,y)=0,
\end{eqnarray*}
where for $i=1,\cdots,m$,
\begin{eqnarray}\label{k}
\kappa_i(x,y,s_i,r,\rho)&:=&z_i(x,y,s_i,r,\rho)+g_i(x,y) +\rho s_i\nonumber\\
&=&\frac{1}{2}[\sqrt{(\rho s_i+ g_i(x,y))^2+4 r \rho}+(\rho s_i+g_i(x,y))].
\end{eqnarray}
We write $z_i(x,y,s_i,r,\rho)$ and $\kappa_i(x,y,s_i,r,\rho)$ as $z_i$ and $\kappa_i$ for convenience, respectively.
Similarly as Liu et al. \cite{ldhs}, properties for $z_i$ and $\kappa_i$ are stated in the following lemma.

\begin{lemma}\label{zkproperty}
For each $r\geq 0, \rho> 0$,  $i=1,\cdots,m$, the following conclusions hold.\\
(1)  $z_i\geq 0$, $\kappa_i\geq 0$, $z_i+g_i(x,y)=\kappa_i-\rho s_i$ and $ z_i \kappa_i=\rho r$;\\
(2)   $g_i(x,y)+z_i=0$ if and only if $g_i(x,y)\leq 0$, $s_i\geq 0$ and $s_i g_i(x,y)=-r$.\\
(3) For $r>0$, $z_i$ and $\kappa_i$ are differentiable with respect to the variable $x, y$, $s$ and $\rho$, respectively,
\begin{eqnarray*}
 \nabla_{(x,y)} z_i=\frac{- z_i}{z_i+\kappa_i} \nabla g_i(x,y),
 && \nabla_{(x,y)} \kappa_i=\frac{ \kappa_i}{z_i+\kappa_i} \nabla g_i(x,y),\\
 \nabla_s z_i=\frac{- \rho z_i}{z_i+\kappa_i}e_i,
 &&\nabla_s \kappa_i=\frac{\rho\kappa_i}{z_i+\kappa_i}e_i.
\end{eqnarray*}
\end{lemma}

\section{Smoothing approximations of the solution mapping}
Let $C^{r,\rho}(x,y,s):=(\phi^{r,\rho}(x,y,s),\psi^{r,\rho}(x,y,s))$.
%
From the Lemma \ref{zkproperty} (1)-(2), for $\rho>0$, the system $C^{r,\rho}(x,y,s)=0$ if and only if 
\begin{eqnarray}
&&\nabla_{y} f(x,y)+  \sum_{i=1}^m s_i \nabla_{y} g_i(x,y) =0,\label{kkt1}\\
&&g_i(x,y)\leq 0,\ s_i\geq 0,\ s_i g_i(x,y)=-r.\label{kkt2}
\end{eqnarray}
For $\bar{r}=0$, any $\bar{\rho}>0$ and a certain point $
\bar x$, if there exists $(\bar y,\bar s)$ such that
$C^{\bar r,\bar \rho}(\bar x,\bar y,\bar s)=0$, then $(\bar y,\bar s)$ is a KKT pair of $(P_{\bar x})$.  
Thus $C^{r,\rho}(x,y,s)=0$ is a perturbation of the KKT system of the problem $(P_{\bar x})$ for any $x$ belongs to a sufficiently small neighborhood of $\bar x$, $r$ near zero and $\rho>0$.
We note that $\bar{\rho}>0$ is an arbitrary constant.

We introduce a family of smoothing functions $\{(y_r^{\rho}(x), s_r^{\rho}(x))\}$, which approximates  $(\bar y,\bar s)$ in Section 3.1 and in Section 3.2, we develop the gradient consistent property, i.e.,
\begin{eqnarray*}
\limsup_{\rho\to\bar{\rho}, r\to 0, x\to \bar x} 
\left (\begin{array}{l}
\nabla y_r^{\rho}(x)\\
\nabla s_r^{\rho}(x)
  \end{array} \right)
  \subseteq 
   M(\bar x).
\end{eqnarray*}


\subsection{Smoothing approximations}

Firstly, it is proven that  $\nabla_{(y,s)}C^{r,\rho}(x,y,s)$ is nonsingular under the Assumption \ref{strongregular}.
\begin{lemma}\label{gamma-inverse}
Assume $(\bar y,\bar u)$ satisfies the Assumption \ref{strongregular} for problem $(P_{\bar x})$.  
Then there exist $\bar{\rho}>0$, $\delta_1>0$, for any $r>0$,  $(r,\rho)\in\mathbb{B}_{\delta_1}(\bar{r},\bar{\rho})$, $\bar r=0$ and $(x,y,s)\in \mathbb{B}_{\delta_1}(\bar x, \bar y,\bar u)$
such that $\nabla_{(y,s)}C^{r,\rho}(x,y,s)$ is nonsingular.
\end{lemma}
{\bf Proof.} Let $w:=(x,y,s)$ and $\bar w:=(\bar x, \bar y,\bar u)$.
Assume to the contrary that there exists a subsequence denoted by $(r,\rho,w)$ converging to $(\bar{r},\bar{\rho},\bar{w})$ such that $\nabla_{(y,s)}C^{r,\rho}(w)$ is singular.
Then there exists a nonzero vector $(\alpha_r^{\rho}(w), \beta_r^{\rho}(w))\in \mathbb{R}^{l+m}$ such that
\begin{eqnarray*}
\nabla_{(y,s)}C^{r,\rho}(w)
 \left (\begin{array}{l}
\alpha_r^{\rho}(w)\\
\beta_r^{\rho}(w)
  \end{array} \right)
  =\left (\begin{array}{ll }
\nabla_y \phi^{r,\rho}(w) & \nabla_s \phi^{r,\rho}(w)\\
\nabla_y \psi^{r,\rho}(w) & \nabla_s \psi^{r,\rho}(w)
  \end{array} \right)
   \left (\begin{array}{l}
\alpha_r^{\rho}(w)\\
\beta_r^{\rho}(w)
  \end{array} \right)
  =0,
\end{eqnarray*}
where
 \begin{eqnarray*}
\nabla_y \phi^{r,\rho}(x,y,s)
&=&\nabla_{yy}^2 f(x,y)+ \sum_{i=1}^m \frac{\kappa_i}{\rho} \nabla_{yy}^2 g_i(x,y)
+\sum_{i=1}^m \frac{\kappa_i}{\rho( z_i+\kappa_i)} \nabla_y g_i(x,y)\nabla_y g_i(x,y)^T,\\
\nabla_s \phi^{r,\rho}(x,y,s)
&=&\sum_{i=1}^m \frac{\kappa_i}{z_i+\kappa_i} \nabla_y g_i(x,y) e_i^T,
\end{eqnarray*}
 \begin{eqnarray*}
\nabla_y \psi^{r,\rho}(x,y,s)
&=& \left (\begin{array}{l}
\frac{\kappa_1}{z_1+\kappa_1} \nabla_y g_1(x,y)^T \\
~~~~~\vdots\\
\frac{\kappa_m}{z_m+\kappa_m} \nabla_y g_m(x,y)^T
  \end{array} \right),\quad
\nabla_s \psi^{r,\rho}(x,y,s)
= \left (\begin{array}{l}
\frac{-\rho z_1}{z_1+\kappa_1}e_1^T  \\
~~~\vdots\\
\frac{-\rho z_m}{z_m+\kappa_m} e_m^T
  \end{array} \right).
\end{eqnarray*}
This equals to
\begin{eqnarray}
&&\nabla_{yy}^2 f(x,y) \alpha_r^{\rho}(w)
+\sum_{i=1}^m \frac{\kappa_i}{\rho( z_i+\kappa_i)} \nabla_y g_i(x,y)\nabla_y g_i(x,y)^T \alpha_r^{\rho}(w) \nonumber\\
&&~~~~~+ \sum_{i=1}^m \frac{\kappa_i}{\rho} \nabla_{yy}^2 g_i(x,y)\alpha_r^{\rho}(w)
+\sum_{i=1}^m \frac{(\beta_r^{\rho}(w))_i\kappa_i}{z_i+\kappa_i} \nabla_y g_i(x,y)=0,\nonumber\\
&& 
\frac{\kappa_i}{z_i+\kappa_i} \nabla_y g_i(x,y)^T \alpha_r^{\rho}(w)
=\frac{\rho (\beta_r^{\rho}(w))_i z_i}{z_i+\kappa_i},\ i=1,\cdots,m,\label{parameter}
\end{eqnarray}
which imply that
\begin{eqnarray}\label{posH}
\nabla_{yy}^2 f(x,y) \alpha_r^{\rho}(w)+ \sum_{i=1}^m \frac{\kappa_i}{\rho} \nabla_{yy}^2 g_i(x,y)\alpha_r^{\rho}(w)
+\sum_{i=1}^m (\beta_r^{\rho}(w))_i \nabla_y g_i(x,y)=0.
\end{eqnarray}

Without loss of generality, suppose there exists $(\alpha,\beta)$ such that 
\begin{eqnarray*}
&&\lim_{(r,\rho,w)\to (\bar{r},\bar{\rho},\bar{w})} (\alpha_r^{\rho}(w), \beta_r^{\rho}(w))/\|(\alpha_r^{\rho}(w), \beta_r^{\rho}(w))\|=(\alpha,\beta).
\end{eqnarray*}
Let $\gamma_{\bar{x}}:=\{i=1,\cdots,m: \bar{u}_i>0\}$.
Since $(\bar y,\bar u)$ is a KKT pair of the problem $(P_{\bar x})$, for the index $i\in \gamma_{\bar{x}}$, $g_i(\bar x, \bar y)=0$, then $z_i\to 0$ and $\kappa_i\to \bar{\rho} \bar{u}_i$ as  $(r,\rho,w)\to (\bar{r},\bar{\rho},\bar{w})$ from the continuity of $g_i$ and the definitions (\ref{z})-(\ref{k}). 
Dividing $\|(\alpha_r^{\rho}(w), \beta_r^{\rho}(w))\|$ on both sides of (\ref{parameter}) and taking limits as $(r,\rho,w)\to (\bar{r},\bar{\rho},\bar{w})$, 
we have that  for $i\in \gamma_{\bar{x}}$,
\begin{eqnarray*}
\nabla_y g_i(\bar x,\bar y)^T \alpha=0.
\end{eqnarray*}
Therefore $\alpha\in {\rm aff} C_{\bar x}(\bar y)$ from the proof of \cite[Lemma 2.2]{dz20}.
Moreover, if $g_i(\bar x,\bar y)<0$, then $z_i>0$ and $\kappa_i\to 0$ as  $(r,\rho,w)\to (\bar{r},\bar{\rho},\bar{w})$ from the continuity of $g_i$ and the definitions (\ref{z})-(\ref{k}). 
Similarly with the above discussion, $\beta_i=0$.

Multiplying $(\alpha_r^{\rho}(w))^T/\|(\alpha_r^{\rho}(w), \beta_r^{\rho}(w))\|^2$ on both sides of (\ref{posH}) and taking limits as $(r,\rho,w)\to (\bar{r},\bar{\rho},\bar{w})$,  then we have that
\begin{eqnarray*}
\alpha^T \nabla_{yy}^2 f(\bar x,\bar y) \alpha+ \sum_{i=1}^m \bar{u}_i \alpha^T \nabla_{yy}^2 g_i(\bar x,\bar y)\alpha
+\sum_{i=1}^m \beta_i \nabla_y g_i(\bar x,\bar y)^T \alpha=0.
\end{eqnarray*}
From (\ref{parameter}), for each $i=1,\cdots,m$, we have that
$
(\beta_r^{\rho}(w))_i \nabla_y g_i(x,y)^T \alpha_r^{\rho}(w)\geq 0
$
and thus for each $i=1,\cdots,m$, $\beta_i \nabla_y g_i(\bar x,\bar y)^T \alpha\geq 0$,
which implies that $\alpha=0$ from the Assumption \ref{strongregular}.

Dividing $\|(\alpha_r^{\rho}(w), \beta_r^{\rho}(w))\|$ on both sides of (\ref{posH}) and taking limits as $(r,\rho,w)\to (\bar{r},\bar{\rho},\bar{w})$, 
$$\sum_{i=1}^m\beta_i \nabla_y g_i(\bar x,\bar y)=0,$$
which implies $\beta=0$ followed from the LICQ in the Assumption \ref{strongregular}.
This is a contraction with the fact that $\|(\alpha,\beta)\|=1$.
Thus there exists a sufficiently small $\delta_1>0$, $\nabla_{(y,s)}C^{r,\rho}(x,y,s)$ is nonsingular
for any $r>0, \rho>0$, $(r,\rho)\in\mathbb{B}_{\delta_1}(\bar{r},\bar{\rho})$  and $(x,y,s)\in \mathbb{B}_{\delta_1}(\bar x, \bar y,\bar u)$.
\BOX

%

We are now ready to investigate the smoothing functions of the primal-dual solution mapping. 
\begin{thm}\label{implicit-c}
Assume $\bar{r}=0$ and $(y(\bar x),u(\bar x))$ satisfies the Assumption \ref{strongregular}  for the problem $(P_{\bar x})$.
Then there exist $\delta>0, \varepsilon>0$, $\bar{\rho}>0$ and Lipschitz continuous mappings
$(y_r^{\rho}(x), s_r^{\rho}(x)):  \mathbb{B}_{\delta}(\bar x)\times \mathbb{B}_{\delta}(\bar{r},\bar{\rho}) \to  \mathbb{B}_{\varepsilon}(y(\bar x),u(\bar x))$ satisfying $C^{r,\rho}(x,y, s)=0$.
Furthermore, for $r>0, \rho>0$, $(y_r^{\rho}(x), s_r^{\rho}(x))$ is continuously differentiable and 
the gradient of $(y_r^{\rho}(x), s_r^{\rho}(x))$ is defined as follows: 
\begin{eqnarray*}
\left (\begin{array}{l}
\nabla y_r^{\rho}(x)\\
\nabla s_r^{\rho}(x)
  \end{array} \right)
 & =&-\nabla_{(y,s)}C^{r,\rho}(x,y_r^{\rho}(x),s_r^{\rho}(x))^{-1} 
\nabla_{x}C^{r,\rho}(x,y_r^{\rho}(x),s_r^{\rho}(x)),
\end{eqnarray*}
where 
$\nabla_{x}C^{r,\rho}(x,y_r^{\rho}(x),s_r^{\rho}(x))=(\nabla_x \phi^{r,\rho}(x,y,s);\nabla_x \psi^{r,\rho}(x,y,s))$.
\end{thm}  
{\bf Proof.} Since 
the system (\ref{kkt1})-(\ref{kkt2}) is strongly regular under the Assumption \ref{strongregular} at $\bar x$, for $r=0$, $\rho>0$,
then the first conclusion of the theorem followed from \cite{r80}. 
We only need to prove the continuously differentiability of $(y_r^{\rho}(x), s_r^{\rho}(x))$ for $r>0, \rho>0$.

From the Lemma \ref{gamma-inverse}, select $\delta<\delta_1$, $\varepsilon<\delta_1$, 
$\nabla_{(y,s)}C^{r,\rho}(x,y,s)$ is nonsingular at $x\in \mathbb{B}_{\delta}(\bar x)$, $(y_r^{\rho}(x), s_r^{\rho}(x))\in \mathbb{B}_{\varepsilon}(y(\bar x), u(\bar x))$ for $r>0, \rho>0$, $(r,\rho)\in\mathbb{B}_{\delta}(\bar{r},\bar{\rho})$.
The differentiability of $(y_r^{\rho}(x),s_r^{\rho}(x))$ follows from \cite[Theorem 9.18]{var}.
Indeed, differentiating both sides of $C^{r,\rho}(x,y_r^{\rho}(x),s_r^{\rho}(x))=0$ with respect to the variable $x$ yields that
\begin{eqnarray}\label{diff-C}
0=\nabla_{x}C^{r,\rho}(x,y_r^{\rho}(x),s_r^{\rho}(x))
+\nabla_{(y,s)}C^{r,\rho}(x,y_r^{\rho}(x),s_r^{\rho}(x))\left (\begin{array}{l}
\nabla y_r^{\rho}(x)\\
\nabla s_r^{\rho}(x)
  \end{array} \right).
\end{eqnarray}
It follows that $(\nabla y_r^{\rho}(x),\nabla s_r^{\rho}(x))^T=-\nabla_{(y,s)}C^{r,\rho}(x,y_r^{\rho}(x),s_r^{\rho}(x))^{-1} 
\nabla_{x}C^{r,\rho}(x,y_r^{\rho}(x),s_r^{\rho}(x))$.
We complete the proof.
\BOX

From easy calculation, $\nabla_x \psi^{r,\rho}(x,y,s)=(\cdots;\frac{\kappa_i}{z_i+\kappa_i} \nabla_x g_i(x,y)^T;\cdots)$ and
$$\displaystyle\nabla_x \phi^{r,\rho}(x,y,s)
=\nabla_{xy}^2 f(x,y)+ \sum_{i=1}^m \frac{\kappa_i}{\rho} \nabla_{xy}^2 g_i(x,y)
+\sum_{i=1}^m \frac{\kappa_i}{\rho( z_i+\kappa_i)} \nabla_y g_i(x,y)\nabla_x g_i(x,y)^T.$$
%
For simplicity, denote by $(\alpha', \beta'):=(\nabla y_r^{\rho}(x),\nabla s_r^{\rho}(x))$.
From the calculations, (\ref{diff-C}) yields that  
\begin{eqnarray}
&&0= \nabla_{xy}^2 f(x,y_r^{\rho}(x)) 
+\sum_{i=1}^m \frac{\kappa_i}{\rho( z_i+\kappa_i)} \nabla_y g_i(x,y_r^{\rho}(x))\nabla_x g_i(x,y_r^{\rho}(x))^T  
+ \sum_{i=1}^m \frac{\kappa_i}{\rho} \nabla_{xy}^2 g_i(x,y_r^{\rho}(x))\nonumber\\
&&~~~+\nabla_{yy}^2 f(x,y_r^{\rho}(x)) \alpha'
+ \sum_{i=1}^m \frac{\kappa_i}{\rho} \nabla_{yy}^2 g_i(x,y_r^{\rho}(x))\alpha'\nonumber\\
&&~~~+\sum_{i=1}^m \frac{\kappa_i}{\rho( z_i+\kappa_i)} \nabla_y g_i(x,y_r^{\rho}(x))\nabla_y g_i(x,y_r^{\rho}(x))^T \alpha'
+\sum_{i=1}^m \frac{\kappa_i }{z_i+\kappa_i} \nabla_y g_i(x,y_r^{\rho}(x))\beta'_i,\nonumber\\
&&0=\frac{\kappa_i}{z_i+\kappa_i} \nabla_y g_i(x,y_r^{\rho}(x))^T \alpha'
-\frac{\rho  z_i}{z_i+\kappa_i}\beta'_i
+\frac{\kappa_i}{z_i+\kappa_i} \nabla_x g_i(x,y_r^{\rho}(x))^T,\ i=1,\cdots,m,\label{nabla2}
\end{eqnarray}
which implies that
\begin{eqnarray}\label{nabla1}
0&=& \nabla_{xy}^2 f(x,y_r^{\rho}(x)) 
+ \sum_{i=1}^m \frac{\kappa_i}{\rho} \nabla_{xy}^2 g_i(x,y_r^{\rho}(x))
+\nabla_{yy}^2 f(x,y_r^{\rho}(x)) \alpha'\nonumber\\
&&+ \sum_{i=1}^m \frac{\kappa_i}{\rho} \nabla_{yy}^2 g_i(x,y_r^{\rho}(x))\alpha'
+\sum_{i=1}^m \nabla_y g_i(x,y_r^{\rho}(x))\beta'_i.
\end{eqnarray}
From (\ref{nabla2}), for each $i=1,\cdots,m$,
\begin{eqnarray*}
\frac{\kappa_i}{\rho} \nabla_y g_i(x,y_r^{\rho}(x))^T \alpha'
-  z_i \beta'_i
+\frac{\kappa_i}{\rho} \nabla_x g_i(x,y_r^{\rho}(x))^T=0,
\end{eqnarray*}
which together with (\ref{nabla1}) implies that $(\alpha',\beta')^T$ is a solution for the system:
$
0=\tilde{b}^{r,\rho}(x,y_r^{\rho}(x),s_r^{\rho}(x))
+\tilde{B}^{r,\rho}(x,y_r^{\rho}(x),s_r^{\rho}(x)) (\alpha',\beta')^T,
$
where 
\begin{eqnarray*}
\tilde{B}^{r,\rho}(x,y,s):=\left (\begin{array}{ll}
\displaystyle\nabla_{yy}^2 f(x,y)+  \sum_{i=1}^m \frac{\kappa_i}{\rho} \nabla_{yy}^2 g_i(x,y)& 
\nabla_{y} g(x,y)^T\\
~~~~~~ (W(x,y,s)-I) \nabla_y g(x,y) & W(x,y,s)
  \end{array} \right)
\end{eqnarray*} 
and
\begin{eqnarray*}
\tilde{b}^{r,\rho}(x,y,s):=\left (\begin{array}{l}
\nabla_{xy}^2 f(x,y)+  \sum_{i=1}^m \frac{\kappa_i}{\rho} \nabla_{xy}^2 g_i(x,y)\\
~~~~~~~~(W(x,y,s)-I) \nabla_x g(x,y)
  \end{array} \right).
\end{eqnarray*} 
Here $W(x,y,s):=diag(w_1(x,y,s),\cdots, w_m(x,y,s))$, $w_i(x,y,s):=\frac{-g_i(x,y)}{s_i-g_i(x,y) }=\frac{ z_i}{ z_i+\kappa_i/\rho}\in [0,1]$, $i=1,\cdots,m$
and $z_i:=z_i(x,y_r^{\rho}(x),(s_r^{\rho}(x))_i, r,\rho)$, $\kappa_i:=\kappa_i(x,y_r^{\rho}(x),(s_r^{\rho}(x))_i, r,\rho)$.
Similarly as the proof of Lemma \ref{gamma-inverse}, $\tilde{B}^{r,\rho}(x,y,s)$ is also nonsingular for any $r>0$,  $(r,\rho)\in\mathbb{B}_{\delta_1}(\bar{r},\bar{\rho})$  and $(x,y,s)\in \mathbb{B}_{\delta_1}(\bar x, \bar y,\bar u)$.

\subsection{The gradient consistent property of the smoothing functions}

From \cite[Proposition 2.3]{dz20}, under Assumption \ref{strongregular}, $\mathcal{A}(x,W)$ is nonsingular, where $W\in \partial_C \Pi_{\mathbb{R}_-^m}(u(x)+g(x,y(x)))$,  for $x\in \mathbb{B}_{\delta_0}(\bar{x})$, $\delta_0$ is defined in the Lemma \ref{implicit-sub}. 
Thus we introduce the following assumption, which is a directly  result under \cite[Proposition 2.3]{dz20} and Lemma \ref{gamma-inverse}.

\begin{ass}\label{C-bound}
We assume that there exist $\delta>0$, $\bar{\rho}>0$, such that
$\|\mathcal{A}(x,W)^{-1}\|$, where $W\in \partial_C \Pi_{\mathbb{R}_-^m}(u(x)+g(x,y(x)))$,
$\|\nabla_{(y,s)}C^{r,\rho}(x,y,s)^{-1}\|$ and $\|\tilde{B}^{r,\rho}(x,y,s)^{-1}\|$
are  bounded with parameter $\mu>0$, for any $r>0$,
$(r,\rho)\in \mathbb{B}_{\delta}(0,\bar{\rho})$ and 
 $(x,y,s)\in \mathbb{B}_{\delta}(\bar x, \bar y,\bar u)$.
\end{ass}

\begin{ass}\label{fass}
We assume the functions $f$ and $g$ satisfy the following properties.
\begin{enumerate}
	\item [(a)]For any $x$ and $y$, the following functions are all Lipschitz continuous with respect to the variables $x$ and $y$, we denote the constants by the same parameter $L$: $\nabla_y f(x,y)$, $\nabla_{xy}^2 f(x,y)$,  $\nabla_{yy}^2 f(x,y)$, 
	$g(x,y)$,
	$\nabla_x g(x,y)$, $\nabla_y g(x,y)$, $\nabla_{xy}^2 g(x,y)$, $\nabla_{yy}^2 g(x,y)$;
	\item [(b)] For any $x$, $y$, there exists $C>0$ such that $\|\nabla_{xy}^2 f(x,y)\|$, $\|g(x,y)\|, \|\nabla_x g(x,y)\|, \|\nabla_y g(x,y)\|$ and $\|\nabla_{xy}^2 g_i(x,y)\|$, $\|\nabla_{yy}^2 g_i(x,y)\|$, $i=1,\cdots,m$ are bounded by $C$.
\end{enumerate}
\end{ass}

In this subsection, we show that the  smoothing function $(y_r^{\rho}(x),s_r^{\rho}(x))$ possesses the gradient consistent property. For a certain point $(\bar x,y(\bar x),u(\bar x))$,
define the index sets:
 \begin{eqnarray*}
 &&\bar{\mathcal{I}}:=\{i: g_i(\bar x,y(\bar x))=0, u_i(\bar x) >0\},\\
 && \bar{\mathcal{J}}:=\{i: g_i(\bar x,y(\bar x))=0, u_i(\bar x) =0\},\\
 && \bar{\mathcal{K}}:=\{i: g_i(\bar x,y(\bar x))<0, u_i(\bar x)=0\}.
\end{eqnarray*}

\begin{thm}\label{distance-v}
Suppose that $(y(\bar x),u(\bar x))$ satisfies the Assumptions \ref{strongregular}, \ref{C-bound}, \ref{fass}.
Then for any $\epsilon>0$, there exists $\delta>0$ such that  for any $r>0, \bar{\rho}>0$, $(r,\rho)\in \mathbb{B}_{\delta}(0,\bar{\rho})$, $x\in \mathbb{B}_{\delta}(\bar x)$,
\begin{eqnarray*}
d\left (
\left (\begin{array}{l}
\nabla y_r^{\rho}(x)\\
\nabla s_r^{\rho}(x)
  \end{array} \right),
M(\bar x)
   \right)
 \leq  \epsilon,
\end{eqnarray*} 
thus
\begin{eqnarray*}
\emptyset \not =\limsup_{\rho\to\bar{\rho}, r\to 0, x\to \bar x} 
\left (\begin{array}{l}
\nabla y_r^{\rho}(x)\\
\nabla s_r^{\rho}(x)
  \end{array} \right)
  \subseteq 
   M(\bar x).
\end{eqnarray*}
If $r=o((g_i(x,y_r^{\rho}(x))+\rho (s_r^{\rho}(x))_i )^2)$ for any $i\in \bar{\mathcal{J}}$,
the set $M(\bar x)$ can be replaced by $M^B(\bar x)$.
\end{thm}
{\bf Proof.}  Assume there exists $c_0>0$ such that $\|(\bar x,y(\bar x),u(\bar x))\|\leq c_0/2$.
Let $v:=(\nabla y_r^{\rho}(x),\nabla s_r^{\rho}(x))^T$, from the Assumptions \ref{C-bound}, \ref{fass},  
there exists $\nu_0>0$ such that
$\|\tilde{b}^{r,\rho}(x,y_r^{\rho}(x),s_r^{\rho}(x))\|\leq \nu_0 C(1+c_0)$ and thus for any $(r,\rho)\in \mathbb{B}_{\delta}(0,\bar{\rho})$, $x\in \mathbb{B}_{\delta}(\bar x)$, $r>0$,
\begin{eqnarray*}
\|v\|\leq \mu \nu_0 C(1+c_0).
\end{eqnarray*}

(i) Without loss of generality, for any $\varepsilon_1>0$, assume there exist $\delta,\varepsilon>0$ and $\widetilde{W}=diag(\widetilde{w}_1,\cdots, \widetilde{w}_m)\in \partial_C \Pi_{\mathbb{R}_-^m}(u(\bar x)+g(\bar x,y(\bar x)))$, where $\widetilde{w}_i=0$ if $i\in \bar{\mathcal{I}}$, $\widetilde{w}_i=1$ if $i\in \bar{\mathcal{K}}$, $\widetilde{w}_i\in [0,1]$, $i\in \bar{\mathcal{J}}$ such that for any $(r,\rho)\in \mathbb{B_{\delta}}(0,\bar{\rho})$ and $(x,y_r^{\rho}(x),s_r^{\rho}(x))\in \mathbb{B_{\delta}}(\bar x)\times\mathbb{B_{\varepsilon}}( y(\bar x),u(\bar x))$, it follows that
$\|W(x,y_r^{\rho}(x),s_r^{\rho}(x))-\widetilde{W}\|\leq \varepsilon_1$ from (\ref{z})-(\ref{k}) and the defintion of $W(x,y,s)$.


Let $\bar{\gamma}:=- \mathcal{A}(\bar x,\widetilde{W})^{-1} a(\bar x,\widetilde{W})\in M(\bar x)$.
Denote $z_i:=z_i(x,y_r^{\rho}(x),(s_r^{\rho}(x))_i, r,\rho)$, $\kappa_i:=\kappa_i(x,y_r^{\rho}(x),(s_r^{\rho}(x))_i, r,\rho)$.
By the Lemma \ref{zkproperty} (1), $\kappa_i/\rho= (s_r^{\rho}(x))_i$. From Assumption \ref{fass}, 
there exists $m_0>0$ such that
\begin{eqnarray*}
&&\|\widetilde{b}^{r,\rho}(x,y_r^{\rho}(x),s_r^{\rho}(x))-a(\bar x,\widetilde{W})\|\\
&\leq& L(2 c_0+5)[\|x-\bar x\|+ \|y_r^{\rho}(x) -y(\bar x)\|]
+2C \|s_r^{\rho}(x) -u(\bar x)\|+ C\varepsilon_1 \\
&\leq& m_0 (\delta + \varepsilon +\varepsilon_1)
\end{eqnarray*} 
and
$
\|\widetilde{B}^{r,\rho}(x,y_r^{\rho}(x),s_r^{\rho}(x))-\mathcal{A}(\bar x,\widetilde{W})\|
\leq m_0 (\delta + \varepsilon +\varepsilon_1).
$
It is easy to see that
\begin{eqnarray*}
0&=&\tilde{b}^{r,\rho}(x,y_r^{\rho}(x),s_r^{\rho}(x))-a(\bar x,\widetilde{W})
+\tilde{B}^{r,\rho}(x,y_r^{\rho}(x),s_r^{\rho}(x))v
-\mathcal{A}(\bar x,\widetilde{W})\bar{\gamma}\\
&=&\tilde{b}^{r,\rho}(x,y_r^{\rho}(x),s_r^{\rho}(x))-a(\bar x,\widetilde{W})
+(\tilde{B}^{r,\rho}(x,y_r^{\rho}(x),s_r^{\rho}(x))
-\mathcal{A}(\bar x,\widetilde{W}))v
+\mathcal{A}(\bar x,\widetilde{W})(v-\bar{\gamma}).
\end{eqnarray*}
Then for any $\epsilon>\mu (1+\mu C\nu_0 (1+c_0))m_0(\delta + \varepsilon +\varepsilon_1)$, 
\begin{eqnarray*} 
\|v-\bar{\gamma}\|&&\leq \|\mathcal{A}(\bar x,\widetilde{W})^{-1}\|
\|\tilde{b}^{r,\rho}(x,y_r^{\rho}(x),s_r^{\rho}(x))-a(\bar x,\widetilde{W})\|\\
&&+ \|\mathcal{A}(\bar x,\widetilde{W})^{-1}\|
\|\tilde{B}^{r,\rho}(x,y_r^{\rho}(x),s_r^{\rho}(x))-\mathcal{A}(\bar x,\widetilde{W})\| \|v\|
\leq  \epsilon.
\end{eqnarray*}

It follows that $(\nabla y_r^{\rho}(x), \nabla s_r^{\rho}(x) )^T \in M(\bar x) +\epsilon  \mathbb{B}(0,1)$, where $\mathbb{B}(0,1)$ is the unit ball centered at $0$.
By the  Theorem \ref{implicit-c}, we have  $\displaystyle\lim_{\rho\to\bar{\rho}, r\to 0, x\to \bar x}  (\nabla y_r^{\rho}(x), \nabla s_r^{\rho}(x) )^T$ exists. 
The compactness of $M(\bar x)$ yields the first conclusion.
%
%
%

(ii) For any $W\in \partial_B \Pi_{\mathbb{R}_-^m}(u(\bar x)+g(\bar x,y(\bar x)))$,  $W=diag(w_1,\cdots, w_m)$, we have that
\begin{eqnarray*}
w_i\in \left \{ \begin{array}{ll}
0,&\mbox{ if } i\in \bar{\mathcal{I}},\\
\{0,1\}, &\mbox{ if }i\in \bar{\mathcal{J}},\\
1,&\mbox{ if } i\in \bar{\mathcal{K}}.
\end{array} \right.
\end{eqnarray*}

For any $i\in \bar{\mathcal{J}}$,  from (\ref{z}) and (\ref{k}), if $r=o((g_i(x,y_r^{\rho}(x))+\rho (s_r^{\rho}(x))_i )^2)$, then $z_i\kappa_i/(z_i+\kappa_i)^2$ converges to zero and thus $(\frac{\kappa_i}{z_i+\kappa_i},\frac{z_i}{z_i+\kappa_i})$ converges to $(0,1)$ or $(1,0)$.
Then the outer limits of $w_i^{r,\rho}(x)$ belong to $\{0,1\}$ and thus we can replace $M(\bar x)$  by $M^B(\bar x)$ in the first conclusion.
\BOX

The following results reveal the distance between  $(y(x),u(x))$ and $(y_r^{\rho}(x),s_r^{\rho}(x))$ for each $r>0$, $\rho>0$.

\begin{lemma}\label{constrainterr}
Assume that the Assumption \ref{fass} holds.
For any $x$, $(y(x),u(x))$ is a KKT pair of $(P_x)$, then for $r\geq 0$, $\rho>0$,
\begin{eqnarray*}\label{bou}
\|C^{r,\rho}(x,y(x),u(x))\|\leq (1+ \frac{ m C}{\rho} ) \sqrt{ r \rho}.
\end{eqnarray*} 
\end{lemma}
 {\bf Proof.}  Let $z^{y,u}:=z(x,y(x),u(x),r,\rho)$, $\kappa^{y,u}:=z(x,y(x),u(x),r,\rho)$.
From the definition (\ref{z}),
\begin{eqnarray*}\label{zi}
z_i^{y,u}= \left \{\begin{array}{ll }
\frac{1}{2}(\sqrt{\rho^2 (u_i(x))^2+4 r \rho}-\rho u_i(x)), & {\rm if}\ u_i(x)>0,\ g_i(x,y(x))=0,\\
\frac{1}{2}(\sqrt{g_i(x,y(x))^2+4 r \rho}-g_i(x,y(x))),& {\rm if}\ u_i(x)=0,\ g_i(x,y(x))<0, \\
\sqrt{ r\rho}, &   {\rm if}\ u_i(x)=0,\ g_i(x,y(x))=0.
  \end{array} \right. 
\end{eqnarray*}
If $u_i(x)>0,\ g_i(x,y(x))=0$,
we have that 
$
\sqrt{\rho^2 (u_i(x))^2+4 r \rho}\leq \rho u_i(x)+2 \sqrt{ r \rho}$
and  if $u_i(x)=0,\ g_i(x,y(x))<0$,
$
\sqrt{g_i(x,y(x))^2+4 r \rho}\leq - g_i(x,y(x))+2 \sqrt{ r\rho}.$
It follows that 
\begin{eqnarray}\label{bou1}
\|z^{y,u}+g(x,y(x))\|\leq  \sqrt{ r \rho}.
\end{eqnarray}
Furthermore, from the KKT condition of $({\rm P}_x)$ and Lemma \ref{zkproperty} (1), we have that
\begin{eqnarray*}
\phi^{r,\rho}(x,y(x),u(x))&=&\nabla_{y} f(x,y(x))+  \sum_{i=1}^m \frac{\kappa_i^{y,u}}{\rho} \nabla_{y} g_i(x,y(x)) \\
&=& \sum_{i=1}^m \frac{1}{\rho} ( \kappa_i^{y,u}-\rho u_i(x) ) \nabla_{y} g_i(x,y(x)) \\
&=& \sum_{i=1}^m \frac{1}{\rho} ( z_i^{y,u}+g_i(x,y(x))) \nabla_{y} g_i(x,y(x)).
\end{eqnarray*}
Thus from the Assumption \ref{fass} (b),
\begin{eqnarray}\label{bou2}
\|\phi^{r,\rho}(x,y(x),u(x))\|\leq
 \frac{ m C}{\rho} \sqrt{ r \rho}.
\end{eqnarray}
Thus the conclusion follows from (\ref{bou1})-(\ref{bou2}) .
%
\BOX

\begin{thm}\label{err}
Assume that the Assumptions \ref{C-bound}, \ref{fass} hold.
Then for any $(r,\rho)\in \mathbb{B}_{\delta}(0,\bar{\rho})$ and $x\in \mathbb{B}_{\delta}(\bar x)$ with $\delta$ and $\bar{\rho}$ are defined in Theorem \ref{implicit-c}, we have that 
\begin{eqnarray*}
\|(y_r^{\rho}(x),s_r^{\rho}(x)) -(y(x),u(x))\| \leq 2\mu  (1+ \frac{ m C}{\rho} ) \sqrt{ r \rho}.
\end{eqnarray*}
\end{thm}  
{\bf Proof.} The conclusion holds automatically if $r=0$.
For $r>0, \rho>0$, we have that
\begin{eqnarray*}
C^{r,\rho}(x,y(x),u(x))&=&C^{r,\rho}(x,y_r^{\rho}(x),s_r^{\rho}(x)) 
+ \nabla_{(y,s)}C^{r,\rho}(x,y_r^{\rho}(x),s_r^{\rho}(x)) (y(x),u(x))- (y_r^{\rho}(x),s_r^{\rho}(x)) \\
&&+o((y(x),u(x))- (y_r^{\rho}(x),s_r^{\rho}(x))),
\end{eqnarray*}
by the Theorem \ref{implicit-c} and the Taylor Expansion.
Since  $C^{r,\rho}(x,y_r^{\rho}(x),s_r^{\rho}(x))=0$, it follows
\begin{eqnarray*}
\|(y(x),u(x))- (y_r^{\rho}(x),s_r^{\rho}(x))\|  \leq 2\mu  (1+ \frac{ m C}{\rho} ) \sqrt{ r \rho}.
\end{eqnarray*}
from Lemma \ref{constrainterr} and the Assumption \ref{C-bound}.
\BOX

\section{Enhanced barrier-smoothing algorithm}
In this section, we replace the solution mapping $y(x)$ in $({\rm SP})$ with its smoothing function $y_r^{\rho}(x)$, which leads to a sequence of approximated single level problems:
\begin{equation}
	\begin{array}{ll}\displaystyle
		\min_{x} &  F(x,y_r^{\rho}(x))\\
		{\rm s.t.} & G(x,y_r^{\rho}(x))\leq 0,\\
		& H(x,y_r^{\rho}(x))= 0.
	\end{array}
	\tag{${\rm SP}_r^{\rho}$}
\end{equation}
We  design a smoothing algorithm which combines the gradient-based method and the augmented Lagrangian method to solve the problem and prove the convergence results.

The augmented Lagrangian function of $({\rm SP}_r^{\rho})$ is defined as
\begin{eqnarray*}
\theta_{\lambda,\mu}^{c}(x,y_r^{\rho}(x))
&:=&F(x,y_r^{\rho}(x))+\frac{1}{2 c}\sum_{i=1}^{p}\left(\max\{0,\lambda_i +c G_i(x,y_r^{\rho}(x))\}^2-\lambda_i^2 \right)\\
&&+\sum_{j=1}^q\left(\mu_j H_j(x,y_r^{\rho}(x))+\frac{c}{2}(H_j(x,y_r^{\rho}(x)))^2\right).
\end{eqnarray*}
Then we consider the unconstrained optimization problem  for $r>0, \rho>0, c>0, \lambda \in \mathbb{R}^{p}, \mu\in \mathbb{R}^{q}$:
\begin{eqnarray*}\label{SCP}
(\rm P_{\lambda,\mu}^c)~~~~~~~\min && \theta_{\lambda,\mu}^{c}(x,y_r^{\rho}(x)). \nonumber
\end{eqnarray*}
In the algorithm, we denote the residual function  measuring the infeasibility and the complementarity by
\begin{eqnarray*}\label{sigma1}
 {\sigma^{\lambda}(x,y_r^{\rho}(x))}  := 
  \max \left \{|H_j(x,y_r^{\rho}(x))|,j=1,\cdots,q, \
|\min\{\lambda_i,-G_i(x,y_r^{\rho}(x))\}|,\ i=1,\cdots,p\right \}.
\end{eqnarray*}
Let $\mu_{min}<0$, $\lambda_{max}, \mu_{max}>0$ be constants and for $\lambda\in \mathbb{R}^p$, $\mu\in \mathbb{R}^q$, $\bar{\lambda}$ and $\bar{\mu}$ be the Euclidean projection of $\lambda, \mu$ onto $\bigotimes_{i=1}^p [0,\lambda_{max}]$ and $\bigotimes_{i=1}^q [\mu_{min},\mu_{max}]$, respectively. 


\begin{alg}\label{algo3-1}
\begin{enumerate} 
\item Given  initial points $(x_1, y_1,s_1)$, 
initial parameters $r_1>0$, $\rho_1>0$, $c_1>0$, initial multipliers
 $\lambda^1,\mu^1$ and tolerance $\varepsilon>0$.
Set \( \{\beta, \delta_1, \delta_2,\bar{\rho}\}\)  within $[0,1)$, $\delta_1<\delta_2$, $\varepsilon_1, \tau_1,  \gamma_1>0$, 
  $\bar{\rho}<\rho_1$.
Set \(k:=1\).

\item If stopping criteria is satisfied, stop the algorithm. 
\item Solve
$\min_{y} f_{r_k}^{\rho_k}(x_k,y,s_k)$
to get an approximate solution $y_{k+1}$ such that
\begin{eqnarray}\label{ycriterion}
\|\nabla_{y} f(x_k,y_{k+1}) + \sum_{i=1}^m \frac{\kappa_i^{k}}{\rho_k}  \nabla_{y} g_i(x_k,y_{k+1}) \|\leq \gamma_k,
\end{eqnarray}
where $z^{k}:=z(x_{k},y_{k+1},s_{k},r_{k},\rho_{k})$ and $\kappa^{k}:=\kappa(x_{k},y_{k+1},s_{k},r_{k},\rho_{k})$.

\item Set
\begin{eqnarray}
 (\tilde{s}_{k+1})_i:=\frac{1}{\rho_k} \kappa_i^{k}=(s_k)_i+ \frac{1}{\rho_k}(z_i^k+g_i(x_k,y_{k+1})),\ i=1,\cdots,m.\label{sdescend}
\end{eqnarray} 
If 
\begin{eqnarray}\label{al4}
\|z(x_{k},y_{k+1},s_{k},r_k,\rho_k)+g(x_{k},y_{k+1})\|\leq \gamma_k,
\end{eqnarray}
set $\gamma_{k+1}:=\delta_1 \gamma_k$, $r_{k+1}:=\delta_1 r_k$, $s_{k+1}:=\tilde{s}_{k+1}$ and
compute $V_k:=(I_l, 0)\tilde{V}_k$, where $\tilde{V}_k$ is computed from
\begin{eqnarray*}\label{inverse}
\tilde{b}^{r_k,\rho_k}(x_k,y_{k+1},s_{k+1})+
\tilde{B}^{r_k,\rho_k}(x_k,y_{k+1},s_{k+1}) \tilde{V}_k=0,
\end{eqnarray*}

Otherwise, set 
 $s_{k+1}:=\tilde{s}_{k+1}$, $(s_{k+1})_i:=-r_k/g_i(x_k,y_{k+1})$ if $g_i(x_k,y_{k+1})<-\varepsilon$, $\rho_{k+1}:= \max\{\bar{\rho}, \delta_2 \rho_k\}$,
 $r_{k+1}:=\delta_2 r_k$,
  $k:=k+1$ and go to the Step 3.

\item  Compute $d_k:=-(\nabla_x \theta_{\bar{\lambda}^k,\bar{\mu}^k}^{c_k}(x_k,y_{k+1})+V_k^T\nabla_y \theta_{\bar{\lambda}^k,\bar{\mu}^k}^{c_k}(x_k,y_{k+1}) )$.
Set $x_{k+1}:=x_k+\alpha_k d_k$, $\tilde{y}_{k+1}=y_{k+1}+\alpha_k V_k d_k$, 
where $\alpha_k:=\beta^{l_k}$, \(l_k\in \{0,1,2\cdots\}\) is the smallest number satisfying
   \begin{eqnarray}
  \theta_{\bar{\lambda}^k,\bar{\mu}^k}^{c_k}(x_{k+1},\tilde{y}_{k+1})- \theta_{\bar{\lambda}^k,\bar{\mu}^k}^{c_k}(x_k,y_{k+1})
   \leq -\alpha_k\delta_{0} \|d_k\|^2.
  \label{ss1}
   \end{eqnarray}

\item  If
    \begin{eqnarray}\label{al3}
\|d_k\|< \tau_k,
    \end{eqnarray}
set
\begin{eqnarray}
    && \lambda_i^{k+1}=\max\{0,\bar{\lambda}_i^k+c_{k} G_i(x_{k},y_{k+1})\},\ i=1,\cdots,p;
   \label{allam}\\
   &&  \mu_{j}^{k+1}=\bar{\mu}_j^k+c_{k} H_j(x_{k},y_{k+1}),\ j=1,\cdots,q\label{almu}
   \end{eqnarray}
and $\tau_{k+1}:=\delta_1 \tau_k$,
 go to the next step.
Otherwise, set $k:=k+1$ and go to Step 3.

\item   If  
  \begin{eqnarray}\label{ale}
   \sigma^{{\lambda}^{k+1}}(x_{k},y_{k+1})<\varepsilon_{k},
   \end{eqnarray}
set $\varepsilon_{k+1}:=\delta_1 \varepsilon_k$, $k:=k+1$ and go to Step 2. 
Otherwise, set $c_{k+1}:=  c_k/\delta_1$, $k:=k+1$ and go to Step 3.
 

\end{enumerate}
\end{alg}

In the algorithm, we first approximate $y_{r_k}^{\rho_k}(x_k)$ by solving the 
 problem $\min_{y} f_{r_k}^{\rho_k}(x_k,y,s_k)$ and then obtain $s_{r_k}^{\rho_k}(x_k)$ through a multiplier update rule when $\gamma_k$ is sufficiently small.
A sufficient decrease in the value of $\theta_{\bar{\lambda},\bar{\mu}}^{c}(x,y_{r}^{\rho}(x))$ is achieved by selecting the stepsize $\alpha_k$, which is based on the replacement of
 $y_{r_k}^{\rho_k}(x_k+\alpha_k d_k)$ with $\tilde{y}_{k+1}=y_{k+1}+\alpha_k V_k d_k$.
This avoids to solve the lower level problem 
 for every $\alpha_k=\beta^{l_k}$, where \(l_k\in \{0,1,2\cdots\}\).
According to Theorem \ref{amijo}, the stepsize $\alpha_k$ exists  for each $k$ and $d_k$ converges to zero as $k$ approaches infinity with reference to Theorem \ref{al-d}, as established by Theorem \ref{al-d}. This ensures that the algorithm can be consistently executed for every iteration.

In order to demonstrate the convergence of Algorithm \ref{algo3-1}, the following assumptions are needed.

\begin{ass}\label{Fass}
We assume the functions $F, G, H$ satisfy the following properties.\\
(a)	$F,G,H$ are twice continuously differentiable;\\
(b)	$\theta_{\bar{\lambda},\bar{\mu}}^{\bar{c}}(x,y)$ is bounded from below, for a fixed $\bar{c}>0$.

\end{ass}



Firstly, we demonstrate that the extended LICQ guarantees that (\ref{al4}) holds for an infinite sequence.
For a point $(\bar x,\bar y)$, the extended LICQ holds at $\bar y$ for $({\rm P}_{\bar x})$ if $\{\nabla_y g_i(\bar x,\bar y): i\in A_{\bar x}(\bar y)\}
$
is linearly independent, where $A_{\bar x}(\bar y):=\{i:g_i(\bar x,\bar y)\geq 0\}$.


\begin{thm}\label{sbound}
Let $\{x_k, y_k, s_k\}$ be the sequence generated by the Algorithm \ref{algo3-1}.
Assume that $(\bar{x},\bar{y})$ is an accumulation point of $\{(x_k, y_{k+1})\}$ and
the extended LICQ holds at $\bar y$ for the problem $({\rm P}_{\bar x})$. 
Then there exists an infinite sequence such that (\ref{al4}) holds and thus $\gamma_k\to 0$ as $k\to \infty$. 
Moreover $\{s_{k+1}\}$ is bounded from above
 for the infinite sequence. 
\end{thm}
{\bf Proof.} Note that from the update rule, $r_k\to 0$ as $k\to \infty$.
 Assume there exists a subsequence $K$ such that $\displaystyle \lim_{k\to\infty, k\in K} (x_k,y_{k+1})=(\bar{x},\bar{y})$. 

(i) Assume to the contrary that there exists $k_0>0$ such that (\ref{al4}) fails for each $k\geq k_0$. 

For the index $i$ such that $g_i(\bar{x},\bar{y})<0$, $(s_{k+1})_i\to 0$ and $(\tilde{s}_{k+1})_i\to 0$ from the update rule of $s_{k+1}$ and the fact that $r_k\to 0$. 
By the definition (\ref{z}), $|z_i^k+g_i(x_k,y_{k+1})|\to 0$ as $k\to \infty$.

For the index $i$ such that $g_i(\bar{x},\bar{y})=0$, we have that $z_i^k\to 0$ and thus $|z_i^k+g_i(x_k,y_{k+1})|\to 0$  as $k\to \infty$ from (\ref{z}) and the fact that $(s_{k})_i\geq 0$.

Then for $k\geq k_0$ large enough, 
there exists at least one index ${i_0}$ and a constant $\beta$ such that $g_{i_0}(\bar{x},\bar{y})>0$ and $z_{i_0}^k+g_{i_0}(x_k,y_{k+1})\geq \beta\gamma_{k_0}$.
From the update rules of $\rho_k$ and $s_k$, for large $k$,
\begin{eqnarray*}
 \frac{1}{\rho_k}(z_{i_0}^k+g_{i_0}(x_k,y_{k+1}))\geq \frac{\beta}{\rho_1}\gamma_{k_0}.
\end{eqnarray*}
Then $(s_{k+1})_{i_0}\to \infty$ as $k\to \infty$. 

Assume without loss of generality there exists $\bar{s}$ such that $\displaystyle \lim_{k\to\infty, k\in K} \frac{s_{k+1}}{\|s_{k+1}\|}=\displaystyle \lim_{k\to\infty, k\in K} \frac{\tilde{s}_{k+1}}{\|s_{k+1}\|}=\bar{s}> 0$ since $\kappa_i^k\to 0$ from (\ref{k}). Dividing by $\|s_{k+1}\|$ on both sides of (\ref{ycriterion}), it follows that
\begin{eqnarray}\label{lower-kkt}
\|\frac{1}{\|s_{k+1}\|} \nabla_{y} f(x_k,y_{k+1}) + \sum_{i=1}^m \frac{\kappa_i^k}{\rho_k\|s_{k+1}\|} \nabla_{y} g_i(x_k,y_{k+1}) \|\leq \ \frac{\gamma_k}{\|s_{k+1}\|}.
 \end{eqnarray}
Taking limits as $k\to\infty, k\in K$, we have that $\bar{s}_i=0$ if $g_i(\bar{x},\bar{y})<0$ and
 \begin{eqnarray}\label{nna1}
 \sum_{i=1}^m \bar{s}_i\nabla_{y} g_i(\bar{x},\bar{y})=0,
 \end{eqnarray}
 which contradicts with the extended LICQ and thus there exists an infinite subset $K_1\subseteq K$ such that (\ref{al4}) holds for each $k\in K_1$, then $\gamma_k\to 0$ from the update rule.
 
 (ii)  Assume to the contrary that $\{s_{k+1}\}_{k\in K_1}$ is unbounded. 
From (\ref{al4}), without loss of generality we assume there exist $\bar{s}$ and $\bar z$
 such that $\displaystyle \lim_{k\to\infty, k\in K_1} \frac{s_{k+1}}{\|s_{k+1}\|}=\bar{s}> 0$, $\lim_{k\to\infty, k\in K_1}z^k=\bar z$ 
and $\bar z+g(\bar x,\bar y)=0$. Hence $g(\bar x,\bar y)\leq 0$.

From (\ref{al4}), for each $i=1,\cdots,m$, $k\in K_1$, $
|z_i^k+ g_i(x_k,y_{k+1}) |\leq  \gamma_k.
$
Multiplying by $\frac{(s_{k+1})_i}{\|s_{k+1}\|}$ on both sides of this inequality,  for $i=1,\cdots,m$, $k\in K_1$, 
\begin{eqnarray*}
\left|\frac{(s_{k+1})_i}{\|s_{k+1}\|}   z_i^k+ \frac{(s_{k+1})_i}{\|s_{k+1}\|} g_i(x_k,y_{k+1}) \right|\leq \frac{(s_{k+1})_i}{\|s_{k+1}\|}\gamma_k.
\end{eqnarray*}
Form Lemma \ref{zkproperty} (1), for each $i=1,\cdots,m$, $\kappa_i^k z_i^k=r_k \rho_k$ and thus by (\ref{sdescend}), $k\in K_1$, 
\begin{eqnarray}\label{sdescend-1}
\left|\frac{(s_{k+1})_i}{\|s_{k+1}\|} g_i(x_k,y_{k+1}) \right|\leq \frac{(s_{k+1})_i}{\|s_{k+1}\|}\gamma_k+\frac{\kappa_i^k}{\rho_k\|s_{k+1}\|}   z_i^k
\leq \gamma_k+\frac{r_k}{\|s_{k+1}\|}.
\end{eqnarray}
Therefore, taking limits as $k\to\infty$, $k\in K_1$ in (\ref{sdescend-1}), we have that
$
\bar{s}_i g_i(\bar x,\bar y)=0.
$

Dividing by $\|s_{k+1}\|$ on both sides of (\ref{ycriterion}), it follows that (\ref{lower-kkt}) holds
and taking limits as $k\to\infty$, $k\in K_1$, we derive that (\ref{nna1}) holds,
  which contradicts with the LICQ and thus $\{s_{k+1}\}_{k\in K_1}$ is bounded.
\BOX

In the rest of this section, we assume that  $K_1$ is the infinite subset such that (\ref{al4}) holds for each $k\in K_1$. 
From the Theorem \ref{sbound}, the sequence $\{s_{k+1}\}_{k\in K_1}$ is bounded under the extended LICQ and the boundedness of $\{(x_k, y_{k+1})\}$, thus in the rest of this section, we assume that $\{(x_k, y_{k+1}, s_{k+1})\}$ is bounded for convenience.
We now investigate the error bound between  $(y_{k+1}, s_{k+1})$ and $(y_r^{\rho}(x_k),s_r^{\rho}(x_k))$ for each $r>0$, $\rho>0$.

\begin{thm}\label{err-solu-r}
Assume that the Assumptions \ref{C-bound}, \ref{fass} hold and $\{(x_k, y_{k+1}, s_{k+1})\}$ is bounded by $M>0$.
Then for $k\in K_1$ large enough,
\begin{eqnarray}\label{serr-1}
\|(y_{k+1},s_{k+1})- (y_{r_k}^{\rho_k}(x_k),s_{r_k}^{\rho_k}(x_k))\|\leq (2\mu +\frac{1}{\bar{\rho}}) \gamma_k.
\end{eqnarray}
Furthermore, 
\begin{eqnarray}\label{serr-2}
\|(y_{k+1},s_{k+1})- (y(x_k),u(x_k))\|\leq (2\mu +\frac{1}{\bar{\rho}}) \gamma_k
+2\mu  (\sqrt{\rho_1}+ \frac{ m C}{\sqrt{\bar{\rho}}} ) \sqrt{ r_k}.
\end{eqnarray}
\end{thm}  
{\bf Proof.}   
For $k\in K_1$, from (\ref{ycriterion}) and (\ref{al4}),
 $\|C^{r_k,\rho_k}(x_k,y_{k+1},s_{k})\|\leq \gamma_k$.
From the Taylor Expansion, for sufficiently large $k\in K_1$,
\begin{eqnarray*}
C^{r_k,\rho_k}(x_k,y_{k+1},s_{k}) 
&=&C^{r_k,\rho_k}(x_k,y_{r_k}^{\rho_k}(x_k),s_{r_k}^{\rho_k}(x_k)) 
 +o((y_{k+1},s_{k})- (y_{r_k}^{\rho_k}(x_k),s_{r_k}^{\rho_k}(x_k)))\\
&&+ \nabla_{(y,s)}C^{r_k,\rho_k}(x_k,y_{r_k}^{\rho_k}(x_k),s_{r_k}^{\rho_k}(x_k)) ( (y_{k+1},s_{k}) - (y_{r_k}^{\rho_k}(x_k),s_{r_k}^{\rho_k}(x_k)) ).
\end{eqnarray*}

From the Assumption \ref{C-bound}, $\|\nabla_{(y,s)}C^{r_k,\rho_k}(x_k, y_{r_k}^{\rho_k}(x_k),s_{r_k}^{\rho_k}(x_k))^{-1}\|$ is bounded by $\mu$ and from the fact that $C^{r_k,\rho_k}(x_k,y_{r_k}^{\rho_k}(x_k),s_{r_k}^{\rho_k}(x_k))=0$, then we have that
\begin{eqnarray}\label{ysdis-0}
\|(y_{k+1},s_{k})- (y_{r_k}^{\rho_k}(x_k),s_{r_k}^{\rho_k}(x_k))\|  \leq  2\mu \gamma_k.
\end{eqnarray}

From (\ref{sdescend}), for sufficiently large $k\in K_1$,
$
\|s_{k+1}-s_k\|\leq \frac{\gamma_k}{\bar{\rho}}.
$
The (\ref{serr-1}) holds by combining this inequality with (\ref{ysdis-0}) and
(\ref{serr-2}) holds from the Theorem \ref{err}.
\BOX

\begin{remark}\label{lipwv}
For each $k$, let $\delta_k:=r_k \rho_k/((g_i(x_k,y_{k+1})+\rho_k (s_{k+1})_i )^2+4r_k \rho_k)$. Assume without loss of generality, there exists $\bar{\delta}\geq 0$ such that $\lim_{k\to\infty}\delta_k=\bar{\delta}$. 
Assume 
$\tilde{V}_k$ generates from the Algorithm \ref{algo3-1}.
Similar to the proof of the Theorem \ref{distance-v}, for any $\epsilon>0$, there exists $\bar{k}>0$ such that for $k\geq \bar{k}, k\in K_1$, we have that $d(\tilde{V}_k,M(x_k))\leq \epsilon$ from (\ref{serr-2}) if $\bar{\delta}>0$ and 
 $d(\tilde{V}_k,M^B(x_k))\leq \epsilon$ if $\bar{\delta}=0$.
\end{remark}

The following theorems establish the existence of the step size such that (\ref{ss1}) holds for each $k$ and demonstrate that $d_k$ converges to zero. Therefore,  it can be concluded that Algorithm \ref{algo3-1} can be effectively executed.
\begin{thm}\label{amijo}
Assume that the Assumptions \ref{fass}, \ref{Fass} hold.
Then for any $k$, there always exists $\alpha_k$ such that (\ref{ss1}) holds.
\end{thm}
{\bf Proof.} Set $x_{k+1}:=x_k+\hat{\alpha} d_k$, $\tilde{y}_{k+1}=y_{k+1}+ \hat{\alpha} V_k d_k$.
For any $k$, from the Taylor Expansion and the definition of $d_k$, since $\theta_{\bar{\lambda}^k,\bar{\mu}^k}^{c_k}(\cdot,\cdot)$ is twice continuously differentiable from the Assumption \ref{Fass},  we have that
\begin{eqnarray*}
\theta_{\bar{\lambda}^k,\bar{\mu}^k}^{c_k}(x_{k+1},\tilde{y}_{k+1})
 &=&\theta_{\bar{\lambda}^k,\bar{\mu}^k}^{c_k}(x_k,y_{k+1})
+\nabla_x \theta_{\bar{\lambda}^k,\bar{\mu}^k}^{c_k}(x_k,y_{k+1})^T (x_{k+1}-x_k)\\
 && +\nabla_y \theta_{\bar{\lambda}^k,\bar{\mu}^k}^{c_k}(x_k,y_{k+1})^T (\tilde{y}_{k+1} - y_{k+1}) +O(\|x_{k+1}-x_k\|^2+ \|\tilde{y}_{k+1} - y_{k+1}\|^2 )\\
 & =&
 \theta_{\bar{\lambda}^k,\bar{\mu}^k}^{c_k}(x_k,y_{k+1})
+\nabla_x \theta_{\bar{\lambda}^k,\bar{\mu}^k}^{c_k}(x_k,y_{k+1})^T \hat{\alpha} d_k\\
 &&
 +[ V_k^T \nabla_y \theta_{\bar{\lambda}^k,\bar{\mu}^k}^{c_k}(x_k,y_{k+1})]^T \hat{\alpha} d_k
+O(\hat{\alpha}^2(1+\|V_k\|^2) \|d_k\|^2 )\\
& =&
 \theta_{\bar{\lambda}^k,\bar{\mu}^k}^{c_k}(x_k,y_{k+1})
-\hat{\alpha}\|d_k\|^2+O(\hat{\alpha}^2(1+\|V_k\|^2) \|d_k\|^2 ).
\end{eqnarray*}
Since $\|V_k\|$ is bounded similar as the proof of the Theorem \ref{distance-v}, then there exists $\alpha_k<\hat{\alpha}$ such that (\ref{ss1}) holds.
\BOX

\begin{thm}\label{al-d} 
Assume that the Assumptions \ref{C-bound}, \ref{fass}, \ref{Fass} hold.
Suppose that the Algorithm \ref{algo3-1} does not terminate within finite iterations.
  Then there exists an infinite subset  such that the condition $(\ref{al3})$ holds.
\end{thm}
{\bf Proof.}  We assume for a contradiction that there exists $\bar k$, for $k\geq \bar{k}$, the condition $(\ref{al3})$ fails and thus there exist $\bar{\varepsilon}>0$, $\bar{\mu}$, $\bar{\lambda}$ and $\bar{c}$ such that for $k> \bar{k}$, $\bar{\mu}^k=\bar{\mu}$, $\bar{\lambda}^k=\bar{\lambda}$, $c_k=\bar{c}$ and
\begin{eqnarray*}
\|d_k\|\geq\bar{\varepsilon}.
\end{eqnarray*}
From the Assumption \ref{Fass}, $\theta_{\bar{\lambda},\bar{\mu}}^{\bar{c}}(x,y)$ is lower bounded and Lipschitz continuous respect to the variable $y$, we assume the  Lipschitz constant is $L_{\theta}$.

From the Theorem \ref{err}, 
$\|y(x_{k+1})-y_{r_k}^{\rho_k}(x_{k+1})\|\leq c_1 \sqrt{ r_k}$, where $c_1=2\mu  (\sqrt{\rho_1}+ \frac{ m C}{\sqrt{\bar{\rho}}} )$.
From the Theorem \ref{implicit-c} and the Taylor Expansion, 
\begin{eqnarray*}
y_{r_k}^{\rho_k}(x_{k+1})=y_{r_k}^{\rho_k}(x_{k}) + \nabla y_{r_k}^{\rho_k}(x_{k}) \alpha_k d_k +o(\|\alpha_k d_k \|).
\end{eqnarray*}
Thus  similar as the proof of the Theorem \ref{distance-v}, from (\ref{serr-1}), there exist $c_2=(2\mu +\frac{1}{\bar{\rho}})$ and $\delta'_k\to 0$ such that
\begin{eqnarray*}
&&\|y(x_{k+1})-\tilde{y}_{k+1}\| \leq 
\|y(x_{k+1})-y_{r_k}^{\rho_k}(x_{k+1})\| + \|y_{r_k}^{\rho_k}(x_{k+1})-\tilde{y}_{k+1}\|\\
&&~~~\leq 
\|y(x_{k+1})-y_{r_k}^{\rho_k}(x_{k+1})\| + \|y_{r_k}^{\rho_k}(x_{k})-y_{k+1}\|
+ \|\nabla y_{r_k}^{\rho_k}(x_{k}) -V_k +o(1 )  \| \| \alpha_k d_k\| \\
&&~~~\leq c_1 \sqrt{ r_k}
+c_2 \gamma_k
+\delta'_k \|d_k\|.
\end{eqnarray*}
Then we have that
\begin{eqnarray}\label{diff-1}
\theta_{\bar{\lambda},\bar{\mu}}^{\bar{c}}(x_{k+1},y(x_{k+1}) )
-\theta_{\bar{\lambda},\bar{\mu}}^{\bar{c}}(x_{k+1},\tilde{y}_{k+1}) )
\leq L_{\theta} [c_1 \sqrt{ r_k}
+c_2 \gamma_k
+\delta'_k \|d_k\|].
\end{eqnarray}
From (\ref{serr-2}),
\begin{eqnarray}\label{diff-2}
\theta_{\bar{\lambda},\bar{\mu}}^{\bar{c}}(x_k,y_{k+1})
-\theta_{\bar{\lambda},\bar{\mu}}^{\bar{c}}(x_{k},y(x_{k}) )
\leq L_{\theta} [c_2 \gamma_k
+c_1 \sqrt{ r_k}].
\end{eqnarray}

Since $r_k, \gamma_k\to 0$ as $k\to \infty$,
 for sufficiently large $k$, from (\ref{ss1}) and (\ref{diff-1})-(\ref{diff-2}), there exists $m_1>0$ such that
\begin{eqnarray*}
\theta_{\bar{\lambda},\bar{\mu}}^{\bar{c}}(x_{k+1},y(x_{k+1}))
-\theta_{\bar{\lambda},\bar{\mu}}^{\bar{c}}(x_{k},y(x_{k}) )
\leq - \alpha_k m_1 \bar{\varepsilon}^2.
\end{eqnarray*}
Since the line search only require a small number of iterations, $\alpha_k$ will never approach to 0, 
which implies that $\theta_{\bar{\lambda},\bar{\mu}}^{\bar{c}}(x_{k+1},y(x_{k+1}))\to -\infty$ as $k\to \infty$.
This contradicts with the boundedness of $\theta_{\bar{\lambda},\bar{\mu}}^{\bar{c}}(x,y)$.
\BOX

We next prove  that any accumulation of the sequence generated by the algorithm is a C-stationary point of (SP), or a B-stationary point under some additional conditions.
Note that if $({\rm P}_{\bar{x}})$ is convex respect to $y$, then the accumulation point becomes to the  C-stationary point or B-stationary point of (BP), respectively.
\begin{thm}\label{alsta}
Assume that the Assumptions \ref{strongregular}, \ref{C-bound}, \ref{fass}, \ref{Fass} hold and $\{(x_k, y_{k+1}, s_{k+1})\}$ is bounded.
Suppose that $\{\lambda_k,\mu_k\}$ is bounded.  
Then the following conclusions hold:\\
(a) Any accumulation point of $\{(x_k, y_{k+1})\}$ is a C-stationary point of (SP).\\
(b) If there exists a sequence $\{\delta_k\}$ converging to zero such that for each $i$ and large $k$, $r_k \rho_k/(g_i(x_k,y_{k+1})+\rho_k (s_{k+1})_i )^2\leq \delta_k$, then any accumulation point of $\{(x_k, y_{k+1})\}$ is a B-stationary point of (SP).
\end{thm}
{\bf Proof.}  Since $\{\lambda_k,\mu_k\}$ is bounded, then $\varepsilon_k\to 0$, which is equivalent to saying that condition $(\ref{ale})$ holds for an infinite subsequence $K_2\subseteq  K_1$ and thus $\displaystyle \lim_{k\rightarrow \infty,k\in K_2} \sigma^{{\lambda}^{k+1}}(x_{k},y_{k+1})= 0$.

Without loss of generality, suppose there exist a subset $K_3\subseteq K_2$, $(\bar{\lambda},\bar{\mu})$ and $(\bar x,\bar y,\bar{s})$ such that $\lim\limits_{k\to\infty, k\in K_3}(x_{k},y_{k+1},s_{k+1})=(\bar x, \bar y,\bar{s})$ and
$\displaystyle\lim_{k\to\infty, k\in K_3}(\lambda_{k+1},\mu_{k+1})=(\bar{\lambda},\bar{\mu})$.
From the definition of $\sigma^{{\lambda}}(x,y)$, it follows that
$H(\bar x,\bar y)=0$, $G(\bar x,\bar y)\leq 0$ and $\lambda_i G_i(\bar x,\bar y)=0$, $i=1,\cdots,p$.


Since $\gamma_k\to 0$, taking limits in (\ref{ycriterion}) and (\ref{al4}),  we have that
\begin{eqnarray*}
&& \nabla_y f(\bar x,\bar y)+   \sum_{i=1}^m \bar{s}_i \nabla_y g_i(\bar x,\bar y)
=0,\label{alkt2}\\
&& \bar{z}+g(\bar x,\bar y)=0,\label{alkt3}
\end{eqnarray*}
where $\bar z$ is the limiting point of $z_k$ without loss of generality
and thus $g_i(\bar x,\bar y)\leq 0$.
Therefore $\bar x$ is a feasible point of (SP) with $y(\bar x)=\bar y$.

From the Remark \ref{lipwv}, without loss of generality assume there exists a vector $\bar{v}$ such that $\lim\limits_{k\to\infty, k\in K_3}V_k=\bar{v}\in M_y(\bar x)$.
From the update rule (\ref{allam})-(\ref{almu}),
we have
\begin{eqnarray}\label{al-kkt}
-d_k&=& \nabla_x \theta_{\bar{\lambda}^k,\bar{\mu}^k}^{c_k}(x_k,y_{k+1})
+V_k^T\nabla_y \theta_{\bar{\lambda}^k,\bar{\mu}^k}^{c_k}(x_k,y_{k+1}) \\
&=&\nabla_x F(x_k,y_{k+1})+\sum_{i=1}^p \lambda_i^{k+1} \nabla_x G_i(x_k,y_{k+1}) 
+\sum_{j=p+1}^q \mu_j^{k+1} \nabla_x H_j(x_k,y_{k+1})\nonumber\\
&+&V_k^T [ \nabla_y F(x_k,y_{k+1})
+\sum_{i=1}^p \lambda_i^{k+1} \nabla_y G_i(x_k,y_{k+1}) 
+\sum_{j=p+1}^q \mu_j^{k+1} \nabla_y H_j(x_k,y_{k+1})].\nonumber
\end{eqnarray}
Since (\ref{al3}) holds, taking  limits  as $k\to\infty$, $k\in K_3$ on (\ref{al-kkt}),   we have that
\begin{eqnarray*}
&&0\in \nabla_x F(\bar x,\bar y)+\sum_{i=1}^p \bar{\lambda}_i \nabla_x G_i(\bar x,\bar y) 
+\sum_{j=p+1}^q \bar{\mu}_j \nabla_x H_j(\bar x,\bar y)\\
&&+M_y(\bar x)^T [\nabla_y F(\bar x,\bar y)
+ \sum_{i=1}^p \bar{\lambda}_i \nabla_y G_i(\bar x,\bar y) 
+\sum_{j=p+1}^q \bar{\mu}_j \nabla_y H_j(\bar x,\bar y)].
\end{eqnarray*}
Then the conclusion (a) holds.

From the Remark \ref{lipwv}, in the case (b), we have that $\lim\limits_{k\to\infty, k\in K_3}V_k\in M_y^B(\bar x)$.
Therefore we derive a B-stationary.
The proof is completed.
\BOX

In the rest of this section, we show the boundedness of the parameters under the ENNAMCQ.
\begin{thm}\label{penaty} 
Assume that $\{(x_k, y_{k+1}, s_{k+1})\}$ is bounded and the Assumptions \ref{strongregular}, \ref{EMF-U}, \ref{C-bound}, \ref{fass}, \ref{Fass} hold 
at any accumulation point of $\{(x_k, y_{k+1})\}$. 
Then $\{ \lambda_k,\mu_k\}$ is bounded and consequently,  the conclusions in Theorem \ref{alsta} hold.
\end{thm}
{\bf Proof.}  
Without loss of generality, suppose there exist a subset $K_2\subseteq K_1$ and $(\bar x,\bar y,\bar{s})$ such that
$\lim\limits_{k\to\infty, k\in K_2}(x_{k},y_{k+1},s_{k+1})=(\bar x, \bar y,\bar{s})$ and (\ref{al3}) holds.

We assume to the contrary that $\|(\lambda_{k},\mu_k)\|\to\infty$ as $k\to \infty$. It follows that $c_k\to \infty$ as $k\to\infty$. There exists a subsequence $\bar{K}_{0}\subseteq K_2$ and $(\lambda,\mu)\in \mathbb{R}^{p+q}$ nonzero such that
$
\lim_{k\to\infty,k\in \bar{K}_{0}}\frac{ (\lambda_{k},\mu_k) }{\|(\lambda_{k},\mu_k)\|}=(\lambda,\mu).
$
It follows from the update rule (\ref{allam}) that $\lambda_i \ge 0,\ i=1,\cdots, p$.

Dividing by $\|(\lambda_{k+1},\mu_{k+1})\|$ on (\ref{al-kkt}) and
taking limits as $k\to \infty, k\in \bar{K}_{0}$,  we have that
\begin{eqnarray}\label{th3.2.0}
&&0\in \sum_{i=1}^p {\lambda}_i \nabla_x G_i(\bar x,\bar y) 
+\sum_{j=p+1}^q {\mu}_j \nabla_x H_j(\bar x,\bar y)\nonumber\\
&&+M_y(\bar x)^T
\left( \sum_{i=1}^p {\lambda}_i \nabla_y G_i(\bar x,\bar y) 
+\sum_{j=p+1}^q {\mu}_j \nabla_y H_j(\bar x,\bar y)\right).
\end{eqnarray}
If $G_i(\bar x,\bar y)<0$, $\lambda_i=0$ from $c_k\to \infty$. 
For each $j=1,\cdots,q$ such that $H_j(\bar x,\bar y)\neq 0$,
 we have for sufficiently large $k\in \bar{K}_{0}$, 
$
\bar{\mu}_j^k H_j(x_{k},y_{k+1})+c_k H_j(x_{k},y_{k+1})^2 >0.
$
Thus
\begin{eqnarray*}
\mu_j H_j(\bar x,\bar y) =\lim_{k\to\infty,k\in \bar{K}_{0}} \frac{\mu_j^{k+1} H_j(x_{k},y_{k+1})}{\|(\lambda_{k+1},\mu_{k+1})\|} 
=\lim_{k\to\infty,k\in \bar{K}_{0}} \frac{\bar{\mu}_j^k H_j(x_{k},y_{k+1})+c_k H_j(x_{k},y_{k+1})^2 }
{\|(\lambda_{k+1},\mu_{k+1})\|}\geq0.
\end{eqnarray*}
Therefore, 
\begin{eqnarray}\label{th3.2.1}
\sum_{i=1}^p {\lambda}_i G_i(\bar x,\bar y)+\sum_{j=1}^q \mu_j H_j(\bar x,\bar y)\geq 0.
\end{eqnarray}
Conditions $(\ref{th3.2.0})$ and $(\ref{th3.2.1})$ contradict with the Assumption \ref{EMF-U}.  From the above discussion and since $(\bar x,\bar y)$ is an arbitrary accumulation point, we know that $\{ \lambda_k,\mu_k\}$ is bounded.
\BOX

\section{Numerical experiments}
To validate the effectiveness of the Enhanced Barrier-Smoothing Algorithm (EBSA), we conducted numerical experiments comparing it with the Levenberg-Marquardt (LM) algorithm proposed by A. Tin and A.B. Zemkoho \cite{tz}, both implemented in MATLAB (R2022b). Numerical tests were based on 134 nonlinear bilevel problems from the current version of the BOLIB library \cite{BOLIB2019}, excluding two problems: one due to the lack of closed-form derivatives and another due to missing extra data. All tests were conducted on a PC with Windows 11 Home system, Intel(R) Core(TM) i9-9900 CPU @ 3.10GHz, and 32 GB of memory.


{\bf Experimental Setup.} In \texttt{EBSA} (Algorithm \ref{algo3-1}), parameters are selected as follows: $\varepsilon=10^{-9},  r_1=1, \rho_1=2, c_1=50, \beta=0.7, \delta_0=0.05, \delta_1=0.8, \delta_2=0.95, \bar{\rho}=10^{-7}, \gamma_1=0.1, \varepsilon_1=0.01, \tau_1=0.8, \lambda^1=\max\{0,c_1G(x_0,y_0)\}, \mu^1=c_1G(x_0,y_0).$

For the LM algorithm, all parameters were chosen according to \cite{tz}. It is worth noting that in our implementation, the LM parameter $\lambda_k$ was updated using the following rule:
\begin{equation*}
	\lambda_k= 0.5\times 1.05^k.
\end{equation*}
The original work \cite{tz} also suggests an alternative approach, which involves selecting the best $\lambda_k$ from a provided set of finite choices. We opted for a fixed update rule to streamline the implementation and reduce computational overhead, as finding suitable candidates for $\lambda_k$ and performing multiple runs can be time-consuming.

{\bf Infeasibility Measure.} 
For a given bilevel problem, let $(\bar{x},\bar{y})$ be the solution produced by a solver. Define the infeasibility measure
\texttt{Infease} as follows:
\begin{equation*}
	\texttt{Infease}:=\| \max(G(\bar{x},\bar{y}),0)\|_\infty
	+ \|H(\bar{x},\bar{y})\|_\infty
	+ \| \max(g(\bar{x},\bar{y}),0)\|_\infty
	+ f(\bar{x},\bar{y})-V(\bar{x}),
\end{equation*}
where $V({x})=\min_y \{f({x},y)\mid g(x,y)\leq 0\}$ is the value function of the lower level problem $(P_x)$, obtained using the MATLAB solver \texttt{fmincon}. \texttt{Infease} is a proper measure of the infeasibility of a bilevel problem at the point $(\bar{x}, \bar{y})$, similarly defined as in \cite{llzz}. A case in the BOLIB library is considered applicable if \texttt{Infease} is smaller than 0.1.

{\bf Initial Points.} 
Initial points were selected as follows.  Let $(\tilde{x}_0,\tilde{y}_0)$ be the initial points suggested by the BOLIB library,
then our initial points $(x_0,y_0)$ were produced by: 
\begin{equation*}
	(x_0,y_0)=(\tilde{x}_0,\tilde{y}_0)+0.01\times(\xi_x,\xi_y),
\end{equation*} 
where $\xi_x$ and $\xi_y$ are random vectors samped from standard normal distributions. Each problem was solved five times starting from independently generated $(x_0,y_0)$ by both \texttt{EBSA} and \texttt{LM}. 

Let  
\begin{equation*}
	R_{F}:=\dfrac{F(\bar{x},\bar{y})-F^*}{1+|F^*|},\quad R_{f}:=\dfrac{f(\bar{x},\bar{y})-f^*}{1+|f^*|},
\end{equation*}
where $F^*$ and $f^*$ are the best known objective function values from the BOLIB library. Missing data are substituted by the smaller objective function values produced by \texttt{EBSA} and \texttt{LM}.

{\bf Stopping Criteria.} 
To decide when to stop \texttt{EBSA}, Let 
\begin{equation*}
	\texttt{Res}_k:=\max\{\|d_k\|,\sigma^{\lambda^k}(x_{k-1},y_{k})\}.
\end{equation*}
When \texttt{Res} is small enough, \texttt{EBSA} should be stopped. In our implementation, \texttt{EBSA} stops when any of the following conditions are met:
\begin{enumerate}
	\item $\texttt{Res}_k<10^{-9}$,
	\item $k>1000$,
	\item $k>200$ and $|\texttt{Res}_k-\texttt{Res}_{k-1}|<10^{-18}$,
	\item $k>300$ and $|\texttt{Res}_k>10^3|$,
	\item $k>300$ and $|\texttt{Res}_k-\texttt{Res}_{k-1}|<10^{-9}$,
	\item $k>800$ and $|\texttt{Res}_k|<10^{-2}$.
\end{enumerate}
For the LM algorithm, the stopping criteria were designed according to \cite{tz}.

{\bf Results and Analysis.} 
Results of the experiments are summarized in Table  \ref{tab:stat}.
The comparison metrics include the number of cases where an algorithm produces smaller $F$ function values (\texttt{\#{}F better}), smaller $f$ function values (\texttt{\#{}f better}), costs less time (\texttt{\#{}time better}), and produces smaller infeasibility (\texttt{\#{}Infease better}). The number of applicable cases is denoted by \texttt{\#{}Applicable cases}. \texttt{avg time} and \texttt{med time} denote the average running time and median running times of the five runs.


 \begin{table}
	\centering
\begin{tabular}{|c|c|c|c|c|c|c|c|}
	\hline
	&  \makebox{\shortstack[c]{\texttt{\#{}Applicable}\\ \texttt{cases}}} &\makebox{\shortstack{\texttt{\#{}F}\\\texttt{better}}} & \makebox{\shortstack{\texttt{\#{}f}\\\texttt{better}}} &\makebox{\shortstack{\texttt{\#{}time}\\\texttt{better}}}  & \makebox{\shortstack{\texttt{\#{}Infease}\\\texttt{better}}}  & \makebox{\shortstack{\texttt{avg}\\\texttt{time}}} & \makebox{\shortstack{\texttt{med}\\\texttt{time}}}  \\
	\hline
	\texttt{EBSA}& 84 &  52  &  48  &  12  &  61  &  1.31e-01  &  7.68e-02    \\
	\hline
	\texttt{LM}& 82 &  30   &  33   &  70   &  20   &  5.44e-02  &  2.98e-02    \\
	\hline
\end{tabular}	
	\caption{Statistical data on the two algorithms\label{tab:stat}}
\end{table}

To better understand the generated data, Figures \ref{fig:Ratios}–\ref{fig:Infeas} are plotted. Ratios $R_F$ and $R_f$, \texttt{Infease}, and \texttt{Time} are sorted in ascending order. 
 Sub-figures on the left show results for all 132 examples, while those on the right focus on the 52 examples where both EBSA and LM found feasible solutions.


\begin{figure}
	\begin{subfigure}{.49\linewidth}
		\pgfplotsset{width=\linewidth}
		\begin{tikzpicture}
			\begin{axis}[
				xlabel={Cases},
				ylabel={$R_F$}, 
				legend entries={\texttt{EBSA},\texttt{LM}},
				legend pos=north west,
				restrict y to domain =-3:4,
				thick
				]					
				\addplot table [x index=0, y index=1, mark=] {fig11.dat};
				\addplot table [x index=0, y index=2, mark=] {fig11.dat};
			\end{axis}
		\end{tikzpicture}
		\caption{All cases}
	\end{subfigure}
	\begin{subfigure}{.49\linewidth}
		\pgfplotsset{width=\linewidth}
		\begin{tikzpicture}
			\begin{axis}[
				xlabel={Cases},
				ylabel={$R_F$}, 
				legend entries={\texttt{EBSA},\texttt{LM}},
				legend pos=north west,
				restrict y to domain =-3:4,
				thick
				]					
				\addplot table [x index=0, y index=1, mark=] {fig12.dat};
				\addplot table [x index=0, y index=2, mark=] {fig12.dat};
			\end{axis}
		\end{tikzpicture}
		\caption{Applicable cases}
	\end{subfigure}
	
	\begin{subfigure}{.49\linewidth}
		\pgfplotsset{width=\linewidth}
		\begin{tikzpicture}
			\begin{axis}[
				xlabel={Cases},
				ylabel={$R_f$}, 
				legend entries={\texttt{EBSA},\texttt{LM}},
				legend pos=north west,
				restrict y to domain =-3:4,
				thick
				]					
				\addplot table [x index=0, y index=1, mark=] {fig13.dat};
				\addplot table [x index=0, y index=2, mark=] {fig13.dat};
			\end{axis}
		\end{tikzpicture}
		\caption{All cases}
	\end{subfigure}
	\begin{subfigure}{.49\linewidth}
		\pgfplotsset{width=\linewidth}
		\begin{tikzpicture}
			\begin{axis}[
				xlabel={Cases},
				ylabel={$R_f$}, 
				legend entries={\texttt{EBSA},\texttt{LM}},
				legend pos=north west,
				restrict y to domain =-3:4,
				thick
				]					
				\addplot table [x index=0, y index=1, mark=] {fig14.dat};
				\addplot table [x index=0, y index=2, mark=] {fig14.dat};
			\end{axis}
		\end{tikzpicture}
		\caption{Applicable cases}
	\end{subfigure}
	\caption{Ratios of objective functions}
	\label{fig:Ratios}
\end{figure}
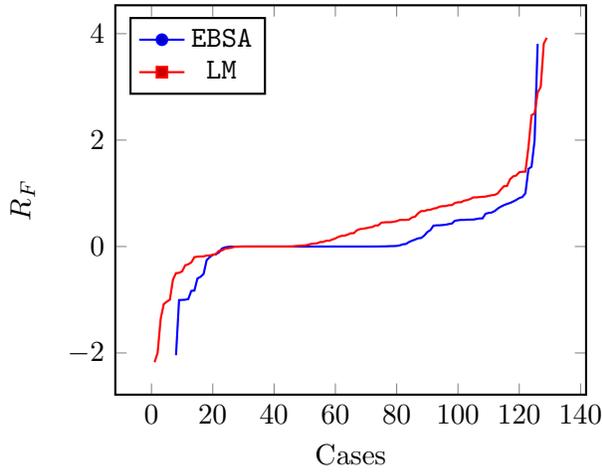
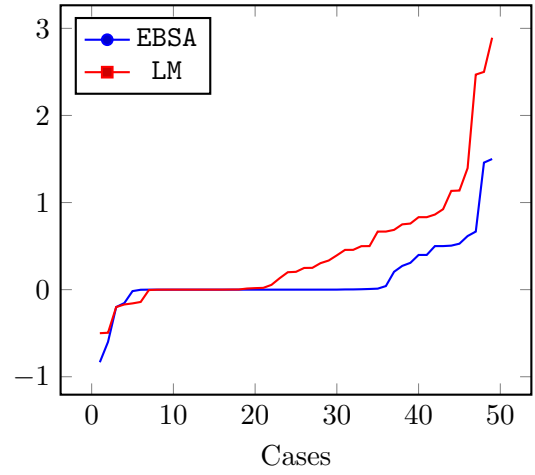
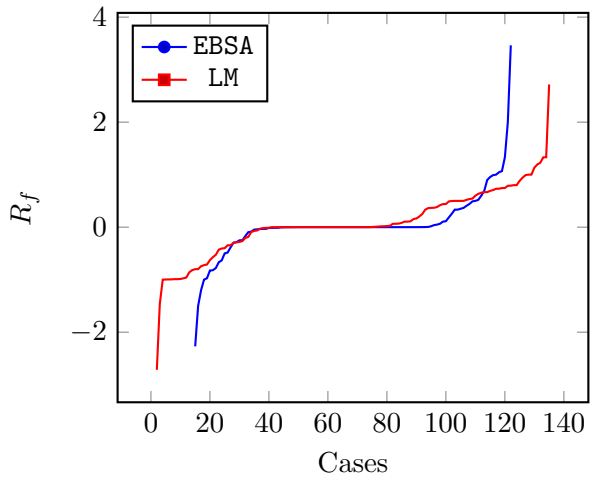
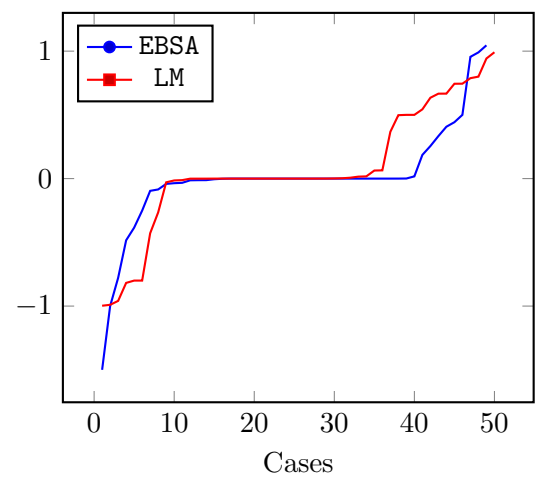

\begin{figure}
	\begin{subfigure}{.49\linewidth}
		\pgfplotsset{width=\linewidth}
		\begin{tikzpicture}
			\begin{semilogyaxis}[
				xlabel={Cases},
				ylabel={\texttt{Time}}, 
				legend entries={\texttt{EBSA},\texttt{LM}},
				legend pos=north west,
				restrict y to domain =-9:1,
				thick
				]					
				\addplot table [x index=0, y index=1, mark=] {fig21.dat};
				\addplot table [x index=0, y index=2, mark=] {fig21.dat};
			\end{semilogyaxis}
		\end{tikzpicture}
		\caption{All cases}
	\end{subfigure}
	\begin{subfigure}{.49\linewidth}
		\pgfplotsset{width=\linewidth}
		\begin{tikzpicture}
			\begin{semilogyaxis}[
				xlabel={Cases},
				ylabel={\texttt{Time}}, 
				legend entries={\texttt{EBSA},\texttt{LM}},
				legend pos=north west,
				restrict y to domain =-9:3,
				thick
				]					
				\addplot table [x index=0, y index=1, mark=] {fig22.dat};
				\addplot table [x index=0, y index=2, mark=] {fig22.dat};
			\end{semilogyaxis}
		\end{tikzpicture}
		\caption{Applicable cases}
	\end{subfigure}

	\caption{Running times on BOLIB}
	\label{fig:Times}
\end{figure}
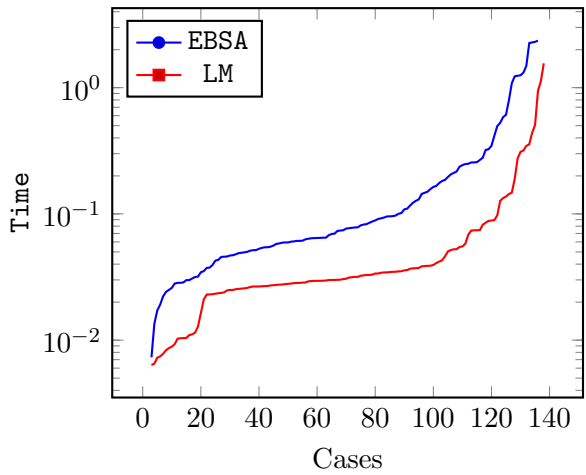
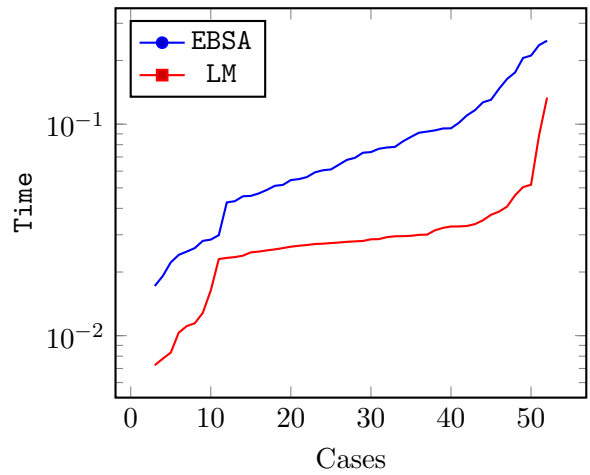

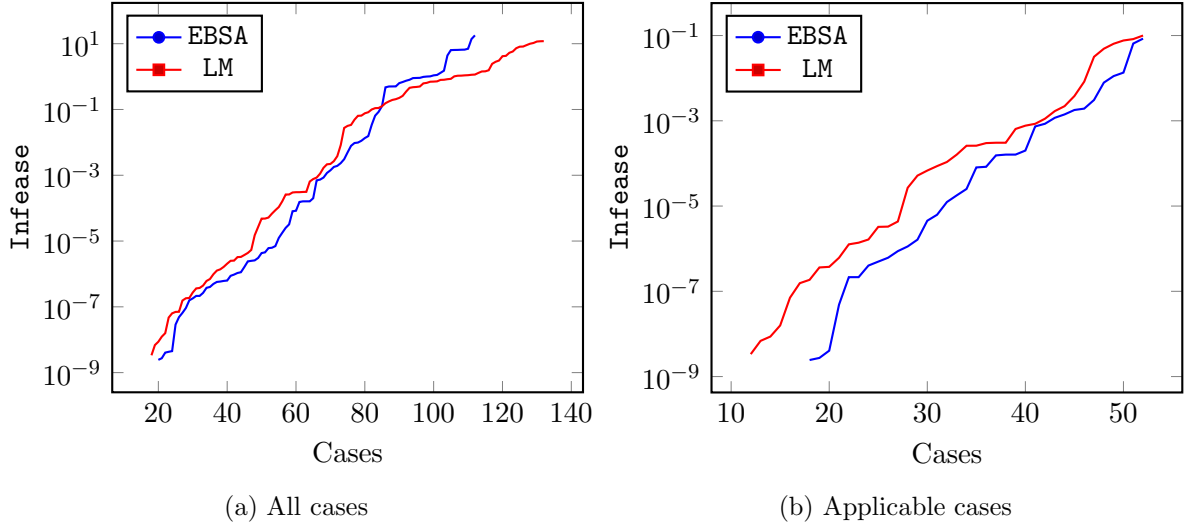
\begin{figure}
	\begin{subfigure}{.49\linewidth}
		\pgfplotsset{width=\linewidth}
		\begin{tikzpicture}
			\begin{semilogyaxis}[
				xlabel={Cases},
				ylabel={\texttt{Infease}}, 
				legend entries={\texttt{EBSA},\texttt{LM}},
				legend pos=north west,
				restrict y to domain =-20:3,
				thick
				]					
				\addplot table [x index=0, y index=1, mark=] {fig31.dat};
				\addplot table [x index=0, y index=2, mark=] {fig31.dat};
			\end{semilogyaxis}
		\end{tikzpicture}
		\caption{All cases}
	\end{subfigure}
	\begin{subfigure}{.49\linewidth}
		\pgfplotsset{width=\linewidth}
		\begin{tikzpicture}
			\begin{semilogyaxis}[
				xlabel={Cases},
				ylabel={\texttt{Infease}}, 
				legend entries={\texttt{EBSA},\texttt{LM}},
				legend pos=north west,
				restrict y to domain =-20:1,
				thick
				]					
				\addplot table [x index=0, y index=1, mark=] {fig32.dat};
				\addplot table [x index=0, y index=2, mark=] {fig32.dat};
			\end{semilogyaxis}
		\end{tikzpicture}
		\caption{Applicable cases}
	\end{subfigure}
	
	\caption{Infeasibility measure}
	\label{fig:Infeas}
\end{figure}

Figures 1-3 provide a detailed comparison:

\begin{enumerate}
	\item {\bf Figure 1.}  Compares upper and lower level objective function values. 
	EBSA consistently achieves lower objective function values at both levels, indicating superior optimization capability.
	
	\item {\bf Figure 2.}  Presents computational time. 
	Although EBSA requires slightly longer computation time, the improvement in objective function values and solution feasibility justifies the additional time.
	
	\item {\bf Figure 3.} Compares infeasibility. 
	EBSA has a lower infeasibility, suggesting higher reliability in finding feasible solutions.

\end{enumerate}

One key advantage of EBSA is that it does not require prior tuning of the parameter, which is often a difficult and time-consuming task. Both algorithms involve solving a system of linear equations, but the adaptive approach to parameter selection in EBSA enhances its practical usability and efficiency in complex problem settings.

Furthermore, the numerical experiments indicate that the LM algorithm demonstrated greater feasibility when the lower level problems were linear-quadratic or nonconvex, as EBSA requires higher convexity in the lower level problem. On the other hand, EBSA demonstrated significant advantages for fractional optimization and higher-dimensional real-world problems.

In conclusion, the numerical experiments confirm the theoretical advantages of EBSA, demonstrating its superior performance over the LM algorithm in terms of feasibility, objective function values, and overall solution quality, particularly for fractional optimization and high-dimensional problems.

%

\section{{Conclusions}}
This paper presented the Enhanced Barrier-Smoothing Algorithm (EBSA) to address the chanllenges  of bilevel optimization problems, particularly those involving nonsmooth solution mappings and hierarchical constraints. By introducing a novel smoothing function for the primal-dual solution mapping, EBSA transforms bilevel problems into a series of smooth single-level problems. This method integrates gradient-based techniques with augmented Lagrangian methods, ensuring convergence to Clarke stationary points and potentially Bouligand stationary points. 
Theoretical analysis and numerical experiments demonstrate the robustness, efficiency, and improved solution accuracy of EBSA. Future works will focus on developing methods that eliminate the need for computing the system of linear equations, further enhancing the efficiency and applicability of the algorithm.


\baselineskip 15pt

\begin{thebibliography} {99}

\bibitem{as} G.B. Allende and G. Still, {\em Solving bilevel programs with the KKT-approach}, Math. Program., {\bf 138}(2013), 309--332.


\bibitem{bkhp} K.P. Bennett, G. Kunapuli, J. Hu and J.-S. Pang, {\em Bilevel optimization and machine learning}, In Computational Intelligence: Research Frontiers: IEEE World Congress on Computational Intelligence, WCCI 2008, Hong Kong, China, June 1-6, 2008, Plenary/Invited Lectures, 25--47, Springer Berlin Heidelberg.










\bibitem{dz20} Y.H. Dai and L. Zhang, {\em Optimality conditions for constrained minimax optimization}, CSIAM Trans. Appl. Math., {\bf 1}(2020), 296--315.






\bibitem{DempeDu} S. Dempe and J. Dutta, {\em Is bilevel programming a special case of a mathematical program with complementarity constraints?} Math. Program., {\bf 131}(2012), 37--48.

\bibitem{df16} S. Dempe and S. Franke, {\em On the solution of convex bilevel optimization problems}, Comput. Optim. Appl., {\bf 63}(2016), 685--703.

\bibitem{dkpk} S. Dempe, V. Kalashnikov, G. P$\acute{e}$rez-Vald$\acute{e}$s and N. Kalashnykova,  {Bilevel Programming Problems}, Energy Systems, Springer Science \& Business Media, Berlin, 2015.

\bibitem{dz} S. Dempe and A.B. Zemkoho, {Bilevel Optimization: Advances and Next Challenges}, Springer, Berlin, 2020.


\bibitem{fzz} A. Fischer, A.B. Zemkoho and S. Zhou, {\em Semismooth Newton-type method for bilevel optimization: Global convergence and extensive numerical experiments}, Optim. Method Softw., {\bf 37}(2022), 1770--1804.

\bibitem{fdfp} L. Franceschi, M. Donini, P. Frasconi and M. Pontil, {\em Forward and reverse gradient-based hyperparameter optimization}, In International Conference on Machine Learning, PMLR, (2017), 1165--1173. 

\bibitem{ffsgp} L. Franceschi, P. Frasconi, S. Salzo, R. Grazzi  and M. Pontil, {\em Bilevel programming for hyperparameter optimization and meta-learning}, In International Conference on Machine Learning, PMLR, (2018), 1568--1577.

\bibitem{gttzz} L.L. Gao, J.J. Ye, H. Yin, S. Zeng and J. Zhang, {\em Moreau envelope based difference-of-weakly-convex reformulation and algorithm for bilevel programs},   arXiv:2306.16761, 2023.




\bibitem{gfps} R. Grazzi, L. Franceschi, M. Pontil and S. Salzo, {\em On the iteration complexity of hypergradient computation}, In International Conference on Machine Learning, PMLR, (2020), 3748--3758.




\bibitem{hss22} E.S. Helou, S.A. Santos and L.E.A. Sim$\tilde{o}$es, {\em A primal nonsmooth reformulation for bilevel optimization problems}, Math.  Program., {\bf 198}(2023), 1381--1409.

\bibitem{hxlt} X. Hu,  N. Xiao, X. Liu and K. C. Toh,   {\em An improved unconstrained approach for bilevel optimization}, SIAM J. Optim., {\bf 33}(2023): 2801--2829.


\bibitem{jmz} L. O. Jolaoso, P. Mehlitz and A. B. Zemkoho, {\em A fresh look at nonsmooth Levenberg–Marquardt methods with applications to bilevel optimization}, Optim., (2024), 1--48.

\bibitem{kpwzs} G. Kornowski, S. Padmanabhan, K. Wang, Z. Zhang and S. Sra, {\em First-Order Methods for Linearly Constrained Bilevel Optimization}, 	arXiv:2406.12771, 2024.


\bibitem{llzz} Y. Li, G.H. Lin, J. Zhang and X. Zhu, {\em A novel approach for bilevel programs based on Wolfe duality}, arXiv:2302.06838, 2023.

\bibitem{llz} Y.-W. Li, G.-H. Lin and X. Zhu, {\em Solving bilevel programs based on lower level Mond–Weir duality}, INFORMS Journal on Computing, (2024). doi:10.1287/ijoc.2023.0108.

\bibitem{lxy} G.-H. Lin, M. Xu and J. J. Ye, {\em On solving simple bilevel programs with a nonconvex lower level program}, Math.  Program., {\bf 144}(2014), 277--305.




\bibitem{lmyzz} R. Liu, P. Mu, X. Yuan, S. Zeng and J. Zhang, {\em A generic first-order algorithmic framework for bi-level programming beyond lower level singleton}, In International Conference on Machine Learning, (2020), 6305--6315.

\bibitem{lmyzz2} R. Liu, P. Mu, X. Yuan, S. Zeng and J. Zhang, {\em A generic descent aggregation framework for gradient-based bi-level optimization}, IEEE Transactions on Pattern Analysis and Machine Intelligence, (2022), 38--57.



\bibitem{ld} X.-W. Liu and Y.-H. Dai, {\em A globally convergent primal-dual interior-point relaxation method for nonlinear programs}, Math. Comp., {\bf 89}(2020), 1301--1329.


\bibitem{ldhs} X.-W. Liu, Y.-H. Dai, Y.-K. Huang and J. Sun, {\em A novel augmented Lagrangian method of multipliers for optimization with general inequality constraints}, Math. Comp., {\bf 92}(2023), 1301--1330.




\bibitem{Liu 24} X. Liu, M. Xu and L.-W. Zhang, {\em Second-order analysis of constrained non-smooth optimization problems with applications to bilevel programs}, submitted. 
    


\bibitem{lvd} J. Lorraine, P. Vicol and D. Duvenaud, {\em Optimizing millions of hyperparameters by implicit differentiation}, In International Conference on Artificial Intelligence and Statistics (AISTATS), (2020), 1540--1552.


\bibitem{lm} Z. Lu and S. Mei,  {\em First-order penalty methods for bilevel optimization},  SIAM J. Optim.,  {\bf 34}(2024), 1937--1969.lm

\bibitem{lpr} Z.-Q. Luo, J.-S. Pang and D. Ralph, {Mathematical Programs with Equilibrium Constraints}, Cambridge University Press, 1996.


\bibitem{Mirrlees99}  J. Mirrlees,  {\em The theory of moral hazard and unobservable behaviour: Part I},  Rev. Econ. Stud., {\bf 66}(1999), 3--22.



\bibitem{nwy} J. Nie, L. Wang and J.J. Ye, {\em Bilevel polynomial programs and semidefinite relaxation methods}, SIAM J. Optim., {\bf 27}(2017), 1728--1757.



\bibitem{Outrata} J.V. Outrata, {\em On the numerical solution of a class of Stackelberg problems}, Z. Oper. Res., {\bf 34}(1990), 255--277.


\bibitem{okz} J. Outrata, M. Kocvara and J. Zowe, {\em Nonsmooth approach to optimization problems with equilibrium constraints: theory, applications and numerical results},  Springer Science \& Business Media, 2013.





\bibitem{rfkl} A. Rajeswaran, C. Finn, S.M. Kakade and S. Levine, {\em Meta-learning with implicit gradients}, in Advances in Neural Information Processing Systems (NeurIPS), (2019),113--124.

\bibitem{r74} S. M. Robinson, {\em Perturbed Kuhn-Tucker points and rates of convergence for a class of nonlinear programming algorithms}, Math. Program., {\bf 7}(1974), 1--16.

\bibitem{r80} S. M. Robinson, {\em Strongly regular generalized equations}, Math. Oper. Res., {\bf 5}(1980), 43--62.

\bibitem{var} R.T. Rockafellar and R.J.-B. Wets, {Variational Analysis}, Springer, Berlin, 1998.

\bibitem{ss17} S. Sabach and S. Shtern, {\em A first order method for solving convex bilevel optimization problems}, SIAM J. Optim., {\bf 27}(2017), 640--660.



\bibitem{s} H. Stackelberg, {Market Structure and Equilibrium}, Springer Science \& Business Media, Berlin, 2010.

\bibitem{tz} A. Tin and A.B. Zemkoho, {\em Levenberg–Marquardt method and partial exact penalty parameter selection in bilevel optimization}, Optim. Eng., {\bf 24}(2023), 1343--1385.


\bibitem{xyz} M. Xu,  J.J. Ye and L. Zhang, {\em Smoothing  sequential quadratic programming method for solving nonconvex, nonsmooth constrained optimization problems}, SIAM J. Optim., {\bf 25}(2015), 1388--1410.



\bibitem{yyzz} J.J. Ye, X. Yuan, S. Zeng and J. Zhang, {\em Difference of convex algorithms for bilevel programs with applications in hyperparameter selection}, Math. Program., {\bf 198}(2023), 1583--1616.



\bibitem{BOLIB2019}
S. Zhou, A.B. Zemkoho and A. Tin, {\em {BOLIB} 2019: bilevel optimization library of test problems version 2}, (2019) arXiv:1812.00230.



\end{thebibliography}
\end{document}